\documentstyle[12pt,righttag]{amsart}
\numberwithin{equation}{section}

\font\germ=eufm10
\newtheorem{thm}{Theorem}[section]
\newtheorem{cor}[thm]{Corollary}
\newtheorem{lem}[thm]{Lemma}
\newtheorem{prop}[thm]{Proposition}

\theoremstyle{definition}
\newtheorem{defn}[thm]{Definition}
\newtheorem{exm}[thm]{Example}

\title[Skew Young diagram method II]
{Skew Young diagram method in spectral decomposition of
integrable lattice models II: Higher levels}

\begin{document}
\newcommand{\ignoreit}[2]{\null}

\maketitle

\begin{center}
{\sc Anatol N.\ Kirillov\footnote[1]{
Permanent address: {\it Steklov Mathematical Institute,
Fontanka 27, St.\ Petersburg, 191011, Russia}},
Atsuo Kuniba$^2$, and Tomoki Nakanishi$^3$}\\
 ~\\
 ~\\
{\it 
$^1$CRM, University of Montreal\\
Succursale A, Centre-ville, Montreal (Quebec) H3W 2P6, Canada\\
 ~\\
$^2$Institute of Physics,
University of Tokyo\\
Komaba, Meguro-ku, Tokyo 153, Japan\\
 ~\\
$^3$Department of Mathematics,
Nagoya University\\
Chikusa-ku, Nagoya 464, Japan}\\
November, 1997
\end{center}

\begin{abstract}
The spectral decomposition of the path space of the vertex model associated
to the level $l$ representation of the quantized affine algebra
$U_q(\widehat{sl}_n)$ is studied.
The spectrum and its degeneracy are parametrized by 
skew Young diagrams and what we call nonmovable tableaux on them,
respectively.
As a result we obtain the characters for the degeneracy of the
spectrum in terms of an alternating sum of skew Schur functions. 
Also studied are 
new combinatorial descriptions (spectral decomposition) of the Kostka numbers
and the Kostka--Foulkes polynomials.
As an application we give a new
proof of Nakayashiki--Yamada's theorem about the branching functions of the
level $l$ basic representation $l\Lambda_k$ of $\widehat{sl}_n$ and 
a generalization of the theorem.
\end{abstract}

\newpage

\section{Introduction}

Study of integrable lattice models often provides
new and remarkable expressions of characters
for various underlying Lie algebras.
For example, the trace of the corner transfer 
matrix \cite{bax} of the {\it
$U_q(\widehat{sl}_n)$-vertex model\/} leads us 
to the formula \cite{DJKMO}
\begin{equation}
\label{eq:intdjkmo}
{\mathrm{ch}}\, {\cal L}(\Lambda_k)(q,x)
=
\sum_{\vec{s}\in {\cal S}_k}
 q^{E(\vec{s})}
e^{{\mathrm{wt}}(\vec{s})},
\qquad k=0,1,\dots,n.
\end{equation}
Here the left hand side is the character of
the level-1 integrable module with highest weight
$\Lambda_k$ of untwisted affine Lie algebra 
$\widehat{sl}_n$.
In the right hand side the summation is taken over the 
{\it spin configurations\/} (={\it path\/}) 
$\vec{s}$; $E(\vec{s})$ and ${\mathrm{wt}}(\vec{s})$
are the energy and the $sl_n$-weight of $\vec{s}$.

There is a natural surjection (the {\it local 
energy map\/}) $\rho : {\cal S}_k
\rightarrow {\mathrm{Sp}}_k$,
where ${\mathrm{Sp}}_k$ denotes the {\it spectrum\/}
of ${\cal S}_k$. The map $\rho$
induces the {\it spectral decomposition\/}
\begin{equation}\label{eq:intch}
{\mathrm{ch}}\, {\cal L}(\Lambda_k)(q,x)
=
\sum_{\vec{h}\in {\mathrm{Sp}}_k}
 q^{E(\vec{h})}\chi_{\vec{h}}(x),
\qquad
\chi_{\vec{h}}(x):=
\sum_{\vec{s}\in \rho^{-1}(\vec{h})}
e^{{\mathrm{wt}}(\vec{s})}.
\end{equation}
This decomposition is remarkable in the
following sense \cite{arakawa,KKN}:
\begin{itemize}
\item[(i)]{
The spectrum ${\mathrm{Sp}}_k$ is parametrized
by the border strips, which are a class 
of skew Young diagrams.
}
\item[(ii)]{
Each subcharacter $\chi_{\vec{h}}$
coincides with the 
skew Schur function associated to the corresponding
border strip.
}
\item[(iii)]{Moreover, $\chi_{\vec{h}}$'s are the
$sl_n$-characters of irreducible
$Y(sl_n)$-modules, where $Y(sl_n)$ is the
 Yangian algebra of type $sl_n$.
}

\end{itemize}

This result admits a simple
representation-theoretical/physical
picture:
The integrals of motions of the hamiltonian
form a commutative algebra ${\cal A}$, and
the algebras ${\cal A}$ and $Y(sl_n)$
act on ${\cal L}(\Lambda_k)$ in a way that
they are mutually commutant.
The formula (\ref{eq:intch}) represents the
reciprocal decomposition of the module ${\cal L}(\Lambda_k)$
with respect to the actions of ${\cal A}$ and $Y(sl_n)$ (cf.\
\cite{T}).

The paper consists of two main parts.

In the first part (Sections \ref{date} and \ref{decom})
we attempt to
generalize the above properties (i)--(iii) to
the levels greater than 1 by considering
the {\it fusion vertex models.}
The result is summarized as follows:
(i) The spectrum is, again, parametrized by
a certain class of skew Young diagrams,
and a formula analogous to (\ref{eq:intch})
is obtained (Theorem \ref{thm:maintheorem}).
(ii) Each subcharacter $\chi_{\vec{h}}$,
also denoted by $t_{\kappa(\vec{h})}$,
is expressed as an alternating sum of skew Schur
functions (Proposition \ref{prop:characterformula}).
(iii)
In the case of either $\widehat{sl}_2$ at higher levels
or $\widehat{sl}_n$ at level 1, those alternating sums
can be reexpressed as a single skew Schur function,
thus reproduces the known results in \cite{arakawa, KKN}.
In the case of $\widehat{sl}_n$ $(n\geq 3)$ at levels greater
than 1, however, they
cannot
be,  in general,
identified with the known characters of irreducible
$Y(sl_n)$-modules, i.e., the characters of the tame modules.
 Thus, the representation-theoretical
interpretation of the decomposition is still missing.

In the second part (Sections \ref{exp} and \ref{kostka})
 we turn into the study of the
{\it Kostka numbers\/} $K_{\lambda,\mu}$
and the {\it Kostka--Foulkes polynomials\/} $K_{\lambda,\mu}(q)$,
\begin{equation*}
K_{\lambda,\mu} =
|{\mathrm{SST}}(\lambda,\mu)|,\qquad
K_{\lambda,\mu}(q)
=\sum_{T\in{\mathrm{SST}}(\lambda,\mu)} q^{c(T)},
\end{equation*}
where ${\mathrm{SST}}(\lambda,\mu)$ is the set of the
semistandard tableaux $T$ of shape $\lambda$
and content $\mu$, and
$c(T)$ is the charge by Lascoux and Sch\"utzenberger.
The set ${\mathrm{SST}}(\lambda,\mu)$ has a natural 
decomposition by the {\it exponents\/} of a tableau,
\begin{equation*}
{\mathrm{SST}}(\lambda,\mu)
=\bigsqcup_{d:{\mathrm{exponents}}}
{\mathrm{SST}}_d(\lambda,\mu).
\end{equation*}
A key observation is the existence of a
 bijection 
$\theta_d : {\mathrm{SST}}_d(\lambda,\mu)
\leftrightarrow
{\mathrm{LR}}_0({\mathrm{Sh}}_d(\mu),\lambda)
$ (Theorem \ref{t5.1}),
where ${\mathrm{LR}}_0({\mathrm{Sh}}_d(\mu),\lambda)$ is another
 set
of the semistandard tableaux satisfying a certain condition.
Thanks to the bijection $\theta_d$, we obtain  new
combinatorial descriptions (spectral decomposition)
of $K_{\lambda,\mu}$
(Corollary \ref{c5.2}) and $K_{\lambda, (l^m)}(q)$
(Theorem \ref{t6.2}), which are the main result in the second
part. 

In Section \ref{sec:truncated} we identify 
the Kostka--Foulkes (and other) polynomials with 
the branching functions of the
{\it truncated characters\/} of the
fusion vertex model.
It follows that the decomposition in
Theorem \ref{t6.2} can be seen as
 the one induced from that in
Theorem \ref{thm:maintheorem}.
As an application we obtain an expression 
of the branching functions of
the integrable $\widehat{sl}_n$-modules as limits of
the Kostka--Foulkes and other polynomials (Corollaries 
\ref{c6.2} and \ref{cor:branchtwo}).

\section{Date-Jimbo-Kuniba-Miwa-Okado correspondence}\label{date}

In this section we review the correspondence between the corner transfer
matrix (CTM) spectra of the vertex model of the symmetric representation 
of $U_q(\widehat{sl}_n)$ and the affine Lie algebra characters of 
$\widehat{sl}_n$ \cite{DJKMO}.

\subsection{Energy function}\label{enf}

Let $\overline{\Lambda}_1,\dots,\overline\Lambda_{n-1}$
 be the fundamental weights
of the Lie algebra $sl_n$, and let
$\epsilon_i=\overline\Lambda_i-\overline\Lambda_{i-1}$
for $i=1,\dots,n$ 
with $\overline\Lambda_0=\overline\Lambda_n=0$.
Then  $\{\epsilon_1,\dots,\epsilon_n\}$ is the set of
the weights of 
the irreducible
representation
 of  $sl_n$ whose highest weight is $\overline\Lambda_1$
(the vector representation).
Let 
$$B_l=\{
v_{a_1\dots a_l} \mid 
1\leq a_1\leq a_2\leq \cdots \leq a_l\leq n\}. 
$$
be a basis of the irreducible representation
 of  $sl_n$ whose highest weight is $l\overline\Lambda_1$
(the $l$-fold symmetric tensor representation) such that 
${\mathrm{wt}}(v_{a_1\dots a_l}) = 
\sum_{i=1}^l \epsilon_{a_i}$.
\par
We define the {\it energy  function}
$H_l: B_l
\times B_l\rightarrow \{0,1,\dots,l\}$ as
\begin{equation}\label{eq:llocalenergy}
H_l(v_{ a_1\dots a_l},
v_{ b_1\dots b_l})=
\min_{\sigma:{\mathrm{permutation}}} \sum_{i=1}^l
H_1(v_{ a_i},v_{ b_{\sigma(i)}}),
\end{equation}
where 
\begin{equation}\label{eq:localenergy}
H_1(v_ a,v_ b)=
\biggl\{
\begin{array}{rl}
0 & \mbox{if $ a <b$}, \\
1 & \mbox{if $a \geq b$}.
\end{array}
\end{equation}
The function $H_l$ is  the logarithm of the $R$-matrix
associated to the $l$-fold symmetric tensor representation
of $U_q(\widehat{sl}_n)$ in the limit $q\rightarrow 0$ \cite{DJM}.

\subsection{Date-Jimbo-Kuniba-Miwa-Okado (DJKMO) correspondence}

For gi\-ven two infinite sequences,
 $\vec{a}=(a_1,a_2,\dots)$
and 
$\vec{b}=(b_1,b_2,\dots)$, of any kind of
objects $a_i,b_i$,
we write $\vec{a}\approx \vec{b}$
if $a_i\neq b_i$ only
 for finitely many $i$.
We often use a shorthand notation
$\vec{a}=(a_1,\dots,a_k,(a_{k+1},\dots,a_{k+m})^\infty)$
for such a periodic sequence as
$\vec{a}=(a_1,\dots,a_k,a_{k+1},
\dots,a_{k+m},a_{k+1},
\dots,a_{k+m},\dots)$.

Let ${\cal K}_l=\{K=((k_1,\ldots ,k_n)^\infty)|
k_i\in {\bold{Z}}_{\geq 0},
\mbox{$\sum_{i=1}^{n}k_{i}=l$. 
}\}$. 
The set ${\cal K}_l$ is identified with the set
of the dominant integral weights of  level $l$
of the untwisted affine Lie algebra
$\widehat{sl}_n$ with the correspondence
\begin{equation*}
K=(k_i)\leftrightarrow \Lambda(K):=k_n\Lambda_0+
\sum_{i=1}^{n-1}k_i\Lambda_i.
\end{equation*}

For $K=(k_i)\in {\cal K}_l$,
we define
$\vec{s}^{(K)}=(s_i^{(K)})_{i=1}^\infty\in B^{\bold{N}}_l$,
where
\begin{equation*}
s_i^{(K)}:=
v_{
\underbrace{\scriptstyle1\dots1}_{\scriptstyle k_{i}}
\underbrace{\scriptstyle2\dots2}_{\scriptstyle k_{i+1}}
\dots
\underbrace{\scriptstyle n\dots n}_{\scriptstyle k_{i+n-1}}},
\quad {\mathrm{wt}}(s^{(K)}_i) = \sum_{j=1}^n k_{i+j-1}\epsilon_j.
\end{equation*}
Then an infinite sequence $\vec{s}=(s_i)_{i=1}^\infty
\in B^{\bold{N}}_l$
 is called a {\it spin configuration} or {\it path} if
$\vec{s}\approx \vec{s}^{(K)}$ for a certain
 $K\in {\cal K}_l$.
The set of the spin configurations 
${\cal S}$, thus, has a natural decomposition
\begin{equation}
{\cal S}=\bigsqcup_{K\in {\cal K}_l} {\cal S}_K,
\quad
{\cal S}_K:=\{ \vec{s}\mid \vec{s}\approx \vec{s}^{(K)}\}.
\end{equation}

For $\vec{s}=(s_i)\in {\cal S}_K$ we define
its {\it energy} $E(\vec{s})$ and $sl_n$-{\it weight}
${\mathrm{wt}}(\vec{s})$ as
\begin{align*}
E(\vec{s})&=\sum_{i=1}^\infty i \left\{H_l(s_{i+1},s_{i})-H_l(s^{(K)}_{i+1},
s^{(K)}_{i})\right\},
\\
{\mathrm{wt}}(\vec{s})&=\sum_{i=1}^{n-1}
k_i\overline\Lambda_{i}
+\sum_{i=1}^\infty \left({\mathrm{wt}}(s_i) - {\mathrm{wt}}(s_i^{(K)})
\right),\quad K=(k_i).
\end{align*}
As is standard in the character theory of $sl_n$,
we regard $e^{{\mathrm{wt}}(\vec{s})}$
as a power of the variables $x_1=e^{\epsilon_1},
x_2=e^{\epsilon_2},\dots,x_n=e^{\epsilon_n}$
with the relation $x_1x_2\cdots x_n=1$.

There is a remarkable connection between
the partition function of ${\cal S}_K$
and an affine Lie algebra character.

\begin{thm}[DJKMO correspondence \cite{DJKMO,kang}] \label{thm:DJKMO}
For a given $K\in {\cal K}_l$,
let ${\cal L}(\Lambda(K))$
be the integrable module of $\widehat{sl}_n$
with the
highest weight $\Lambda(K)$. 
Then the following equality holds:
\begin{equation}\label{eq:djkmo}
{\mathrm{ch}}\, {\cal L}(\Lambda(K))(q,x)
=
\sum_{\vec{s}\in {\cal S}_K}
 q^{E(\vec{s})}
e^{{\mathrm{wt}}(\vec{s})},
\end{equation}
\noindent
where
${\mathrm{ch}}\, {\cal L}(\Lambda(K))$
is the (unnormalized) character of ${\cal L}(\Lambda(K))$
\cite{kac}.
\end{thm}

\subsection{Energy functions and nonmovable tableaux}

Let us describe the energy function $H_l$ in terms of nonmovable tableaux.
For this aim, it is convenient to identify the set
$B_l$ with the crystal of the $l$-fold symmetric tensor
representation of $U_q(sl(n))$.
The latter consists of 
semistandard tableaux
of shape $(l)$ with entries from the set $\{1,2, \ldots, n\}$.
We identify $v_{a_1\ldots a_l}\in B_l$ with a
semistandard tableau as
\begin{equation}
\setlength{\unitlength}{1.2pt}
\begin{picture}(112,10)(-40,-20)
\put(0,-20){\line(1,0){70}}
\put(0,-10){\line(1,0){70}}
\put(0,-10){\line(0,-1){10}}
\put(10,-10){\line(0,-1){10}}
\put(20,-10){\line(0,-1){10}}
\put(60,-10){\line(0,-1){10}}
\put(70,-10){\line(0,-1){10}}
\put(0.5,-20){\makebox(10,10){\small$ a_1$}}
\put(10.5,-20){\makebox(10,10){\small$ a_2$}}
\put(60.5,-20){\makebox(10,10){\small$ a_l$}}
\put(35,-15){\circle*{2}}
\put(40,-15){\circle*{2}}
\put(45,-15){\circle*{2}}
\put(-40,-17){$v_{a_1\dots a_l}=$}
\put(72,-17){.}
\end{picture}
\nonumber
\end{equation}
We shall construct the function
$H_{l_1, l_2}: B_{l_1} \otimes B_{l_2} \rightarrow \mathbf{Z}$
such that $H_l$ in (\ref{eq:llocalenergy}) is realized as 
$H_l = H_{l, l}$ under the above identification.

Next we construct the maps (cf.\ \cite{kang,KN})
$$\widetilde{f}_i: B_l \rightarrow B_l \cup \{ 0 \}, \ \ \widetilde e_i: B_l 
\rightarrow B_l \cup \{ 0 \}, \ \ 1\le i\le n-1.
$$
Let $b\in B_l$ be a semistandard tableau and $i\in\{1,\ldots ,n-1\}$, we 
define $\varphi_i(b)=$ number of $i$ in $b$. For each $i$, $1\le i\le n-1$, 
we define a map
$$\widetilde{f}_i: B_l \rightarrow B_l \cup \{ 0 \}
$$
by the following rule: Let $b'\in B_l$, then 

$\bullet$ $\widetilde f_ib'=0$, if $\varphi_i(b')=0$;

$\bullet$ $\widetilde{f}_i b' = b$ if $\varphi_i(b') > 0$,\\ 
where $b$ is obtained from 
$b'$ by replacing the rightmost $i$ in $b'$ with $i+1$. 

Similarly, we define a map
$$\widetilde e_i:B_l\rightarrow B_l\cup\{ 0\}
$$
by the rule: Let $b\in B_l$, then

$\bullet$ $\widetilde e_ib=0$, if $\varphi_{i+1}(b)=0$;

$\bullet$ $\widetilde e_ib=b'$, if $\varphi_{i+1}(b)>0$,\\
where $b'$ is obtained from $b$ by replacing the leftmost $i+1$ in $b$ 
with $i$.

We also will use the same symbol  $\widetilde{f}_i$ and $\widetilde e_i$, 
$i \in \{1, \ldots , n-1\}$, to denote the maps 
$\widetilde{f}_i: B_{l_1} \otimes B_{l_2}
\rightarrow ( B_{l_1} \otimes B_{l_2}) \cup \{ 0 \}$
and $\widetilde{e}_i: B_{l_1} \otimes B_{l_2}
\rightarrow ( B_{l_1} \otimes B_{l_2}) \cup \{ 0 \}$
which are defined by the rules
\begin{eqnarray}
 \widetilde{f}_i(b_1 \otimes b_2) = \left\{
\begin{array}{ll}
(\widetilde{f}_i b_1 ) \otimes b_2 &
\mbox{ if $\varphi_i(b_1) > \varphi_{i+1}(b_2)$ } \\
b_1 \otimes( \widetilde{f}_i b_2) & \mbox{ otherwise.}
\end{array}\right. \label{1.5}
\end{eqnarray}
\begin{eqnarray}
 \widetilde{e}_i(b_1 \otimes b_2) = \left\{
\begin{array}{ll}
(\widetilde{e}_i b_1 ) \otimes b_2 &
\mbox{ if $\varphi_i(b_1) \ge \varphi_{i+1}(b_2)$ } \\
b_1 \otimes( \widetilde{e}_i b_2) & \mbox{ otherwise.}
\end{array}\right. \label{1.6}
\end{eqnarray}
On the right hand sides of (\ref{1.5}) and (\ref{1.6}) we assume 
$0 \otimes B_{l_2} = B_{l_1} \otimes 0 = 0$.
It is known \cite{kang} that the tensor product of two crystals 
also becomes a crystal under the above rules.

For each $d$, $0 \le d \le \hbox{min}(l_1, l_2)$, we define 
$C_d \subset  B_{l_1} \otimes B_{l_2}$ to be the connected
component of the following (highest weight) element under the actions of 
$\widetilde{f}_i$,  $1\le i\le n-1$,
\begin{equation*}
\begin{picture}(187,25)(0,0)
\put(0,0){\line(1,0){70}}
\put(90,0){\line(1,0){95}}
\put(0,10){\line(1,0){70}}
\put(90,10){\line(1,0){95}}
\put(0,0){\line(0,1){10}}
\put(10,0){\line(0,1){10}}
\put(60,0){\line(0,1){10}}
\put(70,0){\line(0,1){10}}
\put(90,0){\line(0,1){10}}
\put(100,0){\line(0,1){10}}
\put(140,0){\line(0,1){10}}
\put(150,0){\line(0,1){10}}
\put(160,0){\line(0,1){10}}
\put(175,0){\line(0,1){10}}
\put(185,0){\line(0,1){10}}
\put(3,2){\small 1}
\put(30,5){\circle*{2}}
\put(35,5){\circle*{2}}
\put(40,5){\circle*{2}}
\put(63,2){\small 1}
\put(77,2){$\otimes$}
\put(93,2){\small 1}
\put(115,5){\circle*{2}}
\put(120,5){\circle*{2}}
\put(125,5){\circle*{2}}
\put(143,2){\small 1}
\put(153,2){\small 2}
\put(165,5){\circle*{2}}
\put(170,5){\circle*{2}}
\put(178,2){\small 2}
\put(0,11){$\overbrace{\hskip70pt}$}
\put(90,11){$\overbrace{\hskip60pt}$}
\put(151,11){$\overbrace{\hskip35pt}$}
\put(34,19){\small $ l_1$}
\put(95,19){\small $ d+l_2 - l_0$}
\put(157,19){\small $ l_0-d$}
\put(187,0){,}
\end{picture}
\end{equation*}
where $l_0 = \hbox{min}(l_1, l_2)$.
As a set, 
$B_{l_1} \otimes B_{l_2} = 
\bigsqcup_{d=0}^{\min(l_1, l_2)} C_d$.
Moreover as an $U_q(sl(n))$ crystal one has the isomorphism
\begin{eqnarray}
\setlength{\unitlength}{1pt}
\begin{picture}(114,60)(-35,-12)
\put(10,10){\line(1,0){30}}
\put(40,20){\line(1,0){20}}
\put(10,30){\line(1,0){50}}
\put(10,10){\line(0,1){20}}
\put(40,10){\line(0,1){10}}
\put(60,20){\line(0,1){10}}
\put(-35,27){$C_d \simeq B$}
\put(10,30){$\overbrace{\hskip50pt}$}
\put(10,40){\small $\max(l_1,l_2)+d$}
\put(10,10){$\underbrace{\hskip30pt}$}
\put(10,-10){\small $\min(l_1,l_2)-d$}
\put(79,27){,}
\end{picture}
\label{1.7}
\end{eqnarray}
where in the RHS of (\ref{1.7}) the set  
$B_{({\rm max}(l_1,l_2)+d, {\rm min}(l_1,l_2)-d)}$ 
is the $U_q(sl(n))$ crystal described in \cite{KN} 
via the semistandard tableaux
with the specified shape.
We define the energy function 
$H = H_{l_1, l_2}$ : $B_{l_1} \otimes B_{l_2} \rightarrow \mathbf{Z}$
by $H(b_1 \otimes b_2) = d$ for any 
$b_1 \otimes b_2 \in C_d$.
Up to an overall additive constant, this agree with 
the definition from the combinatorial $R$-matrix 
in \cite{NY}.
When $l_1 = l_2$, this definition of the energy function 
coincides with that in (\ref{eq:llocalenergy}) as follows:
\begin{equation*}
H_{l, l}(v_{a_1 \ldots a_l}\otimes v_{b_1 \ldots b_l})
= H_{l}(v_{a_1 \ldots a_l}, v_{b_1 \ldots b_l}).
\end{equation*}

Now we are ready to give a description of the energy function 
$H:=H_{l_1,l_2}$ in terms of nonmovable tableaux.

\begin{defn}\label{d1.1} A skew shape tableau $T$ is called to be 
{\it nonmovable}, if the following conditions are satisfied:

$\bullet$ $T$ is semistandard.

$\bullet$ For each $i$, let us denote by $T_i$ the tableau which is 
obtained from $T$ by moving the $i$th row of $T$ from the right to 
the left by one box without changing the positions of all other rows. 
Then the tableau $T_i$ is either of non-skew shape or non-semistandard.
\end{defn}
 
We denote by 
${\mathrm{NMT}}(\lambda)$ 
the set of all
the nonmovable tableaux of shape $\lambda$ with entries not exceeding $n$.

To continue, let us introduce one more function on the same space
$D = D_{l_1, l_2}$~: $B_{l_1} \otimes B_{l_2} \rightarrow \mathbf{Z}$.
Let $b_1 \in B_{l_1}$ and $b_2 \in B_{l_2}$ be  semistandard 
tableaux.
Then $D = D(b_1\otimes b_2)$ is 
uniquely determined by the condition that 
the following tableau is nonmovable:
\begin{equation*}
\begin{picture}(70,25)(0,1)
\setlength{\unitlength}{1.1pt}
\put(0,0){\line(1,0){40}}
\put(0,10){\line(1,0){70}}
\put(20,20){\line(1,0){50}}
\put(0,0){\line(0,1){10}}
\put(20,10){\line(0,1){10}}
\put(70,10){\line(0,1){10}}
\put(40,0){\line(0,1){10}}
\put(16,2){\small $b_2$}
\put(40,12){\small $b_1$}
\put(0,10){$\overbrace{\hskip21pt}$}
\put(5,18){\small $D$}
\end{picture}
\qquad \mbox{for $l_1 \ge l_2$},
\qquad\qquad
\begin{picture}(70,25)(0,1)
\setlength{\unitlength}{1.1pt}
\put(0,0){\line(1,0){50}}
\put(0,10){\line(1,0){70}}
\put(30,20){\line(1,0){40}}
\put(0,0){\line(0,1){10}}
\put(30,10){\line(0,1){10}}
\put(70,10){\line(0,1){10}}
\put(50,0){\line(0,1){10}}
\put(22,2){\small $b_2$}
\put(48,12){\small $b_1$}
\put(51,10){$\underbrace{\hskip21pt}$}
\put(55,-3){\small $D$}
\end{picture}
\qquad \mbox{for $l_1 \le l_2$}.
\end{equation*}
It is clear that $0 \le D_{l_1, l_2} \le {\rm min }(l_1,l_2)$.

\begin{prop}\label{p1.1}
\begin{equation}
H_{l_1,l_2} = D_{l_1,l_2}. \label{1.8}
\end{equation}
\end{prop}
\begin{pf}
We assume $l_1 \ge l_2$.
The other case is similar.
First we check the statement for the highest weight elements.
In the component $B_{(l_1 + d, l_2-d)}$ it is given by
\begin{equation*}
\begin{picture}(162,25)(0,0)
\put(0,0){\line(1,0){70}}
\put(90,0){\line(1,0){70}}
\put(0,10){\line(1,0){70}}
\put(90,10){\line(1,0){70}}
\put(0,0){\line(0,1){10}}
\put(10,0){\line(0,1){10}}
\put(60,0){\line(0,1){10}}
\put(70,0){\line(0,1){10}}
\put(90,0){\line(0,1){10}}
\put(100,0){\line(0,1){10}}
\put(115,0){\line(0,1){10}}
\put(125,0){\line(0,1){10}}
\put(135,0){\line(0,1){10}}
\put(150,0){\line(0,1){10}}
\put(160,0){\line(0,1){10}}
\put(3,2){\small 1}
\put(30,5){\circle*{2}}
\put(35,5){\circle*{2}}
\put(40,5){\circle*{2}}
\put(63,2){\small 1}
\put(77,2){$\otimes$}
\put(93,2){\small 1}
\put(105,5){\circle*{2}}
\put(110,5){\circle*{2}}
\put(118,2){\small 1}
\put(128,2){\small 2}
\put(140,5){\circle*{2}}
\put(145,5){\circle*{2}}
\put(153,2){\small 2}
\put(0,11){$\overbrace{\hskip70pt}$}
\put(90,11){$\overbrace{\hskip35pt}$}
\put(125,11){$\overbrace{\hskip35pt}$}
\put(34,19){\small $l_1$}
\put(106,19){\small $d$}
\put(132,19){\small $l_2-d$}
\put(162,0){.}
\end{picture}
\end{equation*}
Correspondingly, the tableau 
\begin{equation*}
\begin{picture}(105,35)(0,0)
\put(0,0){\line(1,0){70}}
\put(0,10){\line(1,0){105}}
\put(35,20){\line(1,0){70}}
\put(0,0){\line(0,1){10}}
\put(10,0){\line(0,1){10}}
\put(25,0){\line(0,1){10}}
\put(35,0){\line(0,1){20}}
\put(45,0){\line(0,1){20}}
\put(60,0){\line(0,1){10}}
\put(70,0){\line(0,1){10}}
\put(95,10){\line(0,1){10}}
\put(105,10){\line(0,1){10}}
\put(3,2){\small 1}
\put(28,2){\small 1}
\put(38,2){\small 2}
\put(63,2){\small 2}
\put(38,12){\small 1}
\put(98,12){\small 1}
\put(15,5){\circle*{2}}
\put(20,5){\circle*{2}}
\put(50,5){\circle*{2}}
\put(55,5){\circle*{2}}
\put(65,15){\circle*{2}}
\put(70,15){\circle*{2}}
\put(75,15){\circle*{2}}
\put(0,11){$\overbrace{\hskip34pt}$}
\put(35,21){$\overbrace{\hskip69pt}$}
\put(16,19){\small $d$}
\put(69,29){\small $l_1$}
\end{picture}
\end{equation*}
is indeed nonmovable.
Now we have 
only to show $D(\widetilde{f}_i(b_1 \otimes b_2)) = D(b_1 \otimes b_2)$
for $\widetilde{f}_i(b_1 \otimes b_2) \neq 0$.
Put $s = \varphi_i(b_1)$ and 
$s' = \varphi_{i+1}(b_2)$.
We consider the cases (i) $s > s'$ and  (ii) $s \le s'$ separately.
In the case (i), $\widetilde{f}_i(b_1 \otimes b_2) = (\widetilde{f}_i b_1)
\otimes b_2$
and $\widetilde{f}_i b_1$ is obtained by changing the rightmost $i$
in $b_1$ into $i+1$. 
Thus $D(\widetilde{f}_i(b_1 \otimes b_2)) \ge D(b_1 \otimes b_2)$. 
Suppose $D(\widetilde{f}_i(b_1 \otimes b_2)) > D(b_1 \otimes b_2)$. 
This can happen only in the situation
\begin{equation*}
\begin{picture}(157,52)(0,-15)
\put(0,0){\line(1,0){130}}
\put(0,10){\line(1,0){155}}
\put(20,20){\line(1,0){135}}
\put(0,0){\line(0,1){10}}
\put(20,10){\line(0,1){10}}
\put(25,10){\line(0,1){10}}
\put(30,0){\line(0,1){10}}
\put(50,10){\line(0,1){10}}
\put(55,0){\line(0,1){10}}
\put(75,0){\line(0,1){20}}
\put(100,0){\line(0,1){20}}
\put(125,0){\line(0,1){20}}
\put(130,0){\line(0,1){10}}
\put(155,10){\line(0,1){10}}
\put(32,2){\small $i+1$}
\put(60,5){\circle*{2}}
\put(65,5){\circle*{2}}
\put(70,5){\circle*{2}}
\put(77,2){\small $i+1$}
\put(111,2){\small $r'$}
\put(37,12){\small $i$}
\put(55,15){\circle*{2}}
\put(60,15){\circle*{2}}
\put(65,15){\circle*{2}}
\put(70,15){\circle*{2}}
\put(87,12){\small $i$}
\put(112,12){\small $r$}
\put(93,18){\vector(3,2){15}}
\put(107,30){\small $i+1$}
\put(30,-1){$\underbrace{\hskip67pt}$}
\put(60,-15){\small $s'$}
\put(25,21){$\overbrace{\hskip72pt}$}
\put(59,29){\small $s$}
\put(157,0){,}
\end{picture}
\end{equation*}
where $i+1 \le r < r'$.
However from $s > s'$ the tableau here for $b_1 \otimes b_2$ 
is already non-semistandard, which is a contradiction.
In the case (ii), $\widetilde{f}_i(b_1 \otimes b_2) =  b_1 \otimes 
(\widetilde{f}_i b_2)$
and $\widetilde{f}_i b_2$ is obtained by changing the rightmost $i$
in $b_2$ into $i+1$.
Thus $D(\widetilde{f}_i(b_1 \otimes b_2)) \le D(b_1 \otimes b_2)$. 
Suppose $D(\widetilde{f}_i(b_1 \otimes b_2)) <  D(b_1 \otimes b_2)$. 
This can happen only in the situation
\begin{equation*}
\begin{picture}(177,51)(0,-16)
\put(0,0){\line(1,0){150}}
\put(0,10){\line(1,0){175}}
\put(20,20){\line(1,0){155}}
\put(0,0){\line(0,1){10}}
\put(20,10){\line(0,1){10}}
\put(25,0){\line(0,1){20}}
\put(50,0){\line(0,1){20}}
\put(75,0){\line(0,1){20}}
\put(95,10){\line(0,1){10}}
\put(100,0){\line(0,1){10}}
\put(125,0){\line(0,1){10}}
\put(120,10){\line(0,1){10}}
\put(145,10){\line(0,1){10}}
\put(150,0){\line(0,1){10}}
\put(175,10){\line(0,1){10}}
\put(37,12){\small $r$}
\put(62,12){\small $i$}
\put(80,15){\circle*{2}}
\put(85,15){\circle*{2}}
\put(90,15){\circle*{2}}
\put(107,12){\small $i$}
\put(132,12){\small $r'$}
\put(37,2){\small $i$}
\put(52,2){\small $i+1$}
\put(82,5){\circle*{2}}
\put(87,5){\circle*{2}}
\put(92,5){\circle*{2}}
\put(102,2){\small $i+1$}
\put(34,2){\vector(-3,-2){15}}
\put(10,-16){\small $i+1$}
\put(50,21){$\overbrace{\hskip70pt}$}
\put(84,29){\small $s$}
\put(50,-1){$\underbrace{\hskip75pt}$}
\put(84,-15){\small $s'$}
\put(177,0){,}
\end{picture}
\end{equation*}
where $r < i < r'$. 
In fact 
only $s = s'$ is possible so that 
the above tableau for $b_1 \otimes b_2$ 
becomes semistandard within the constraint $s \le s'$.
But then $D(\widetilde{f}_i(b_1 \otimes b_2)) <  D(b_1 \otimes b_2)$
implies $r' < i+1$, which is a contradiction.
\end{pf}

Let us give another and
a more elementary proof of Proposition~\ref{p1.1} 
in the homogeneous case $l_1=l_2=l$ using the definition of
(\ref{eq:llocalenergy}).
For any pair $v_{ a_1\dots a_l},
v_{ b_1\dots b_l}\in B_l$ and 
an integer $d\in \{0,1,\dots,l\}$, let
$T_d(v_{ a_1\dots a_l},v_{ b_1\dots b_l})$
be the tableau below:
\begin{equation}\label{eq:tskew}
\setlength{\unitlength}{1.2pt}
\begin{picture}(186,48)(-84,-33)
\put(0,-20){\line(1,0){70}}
\put(0,-10){\line(1,0){100}}
\put(30,0){\line(1,0){70}}
\put(0,-10){\line(0,-1){10}}
\put(10,-10){\line(0,-1){10}}
\put(20,-10){\line(0,-1){10}}
\put(30,0){\line(0,-1){20}}
\put(40,0){\line(0,-1){12}}
\put(40,-18){\line(0,-1){2}}
\put(50,0){\line(0,-1){10}}
\put(60,0){\line(0,-1){20}}
\put(70,0){\line(0,-1){2}}
\put(70,-8){\line(0,-1){12}}
\put(90,0){\line(0,-1){10}}
\put(100,0){\line(0,-1){10}}
\put(0.5,-20){\makebox(10,10){\small$ b_1$}}
\put(10.5,-20){\makebox(10,10){\small$ b_2$}}
\put(60.5,-20){\makebox(10,10){\small$ b_l$}}
\put(30.5,-20){\makebox(20,10){\small$ b_{d+1}$}}
\put(60.5,-20){\makebox(10,10){\small$ b_l$}}
\put(30.5,-10){\makebox(10,10){\small$ a_1$}}
\put(40.5,-10){\makebox(10,10){\small$ a_2$}}
\put(60.5,-10){\makebox(20,10){\small$ a_{l-d}$}}
\put(90.5,-10){\makebox(10,10){\small$ a_l$}}
%
\put(0,-9){$\overbrace{\hskip33pt}$}
\put(11,-1){\small $d$}
\put(30,1){$\overbrace{\hskip83pt}$}
\put(63,9){\small $l$}
\put(0,-20){$\underbrace{\hskip83pt}$}
\put(33,-33){\small $l$}
\put(-84,-10){$T_d(v_{ a_1\dots a_l},v_{ b_1\dots b_l})=$}
\put(102,-10){.}
\end{picture}
\end{equation}

\begin{prop}\label{prop:hbytableau}
Let  $v_{ a_1\dots a_l},
v_{ b_1\dots b_l}\in B_l$ and $d\in\{0,1,\dots,l\}$. Then 
\hfill\break\noindent 
(i) $H_l(v_{ a_1\dots a_l},v_{ b_1\dots b_l})\leq d$
if and only if the tableau $T_d(v_{ a_1\dots a_l},
v_{ b_1\dots b_l})$ is semistandard.
\hfill\break\noindent 
(ii) $H_l(v_{ a_1\dots a_l},v_{ b_1\dots b_l}) =  d$
if and only if the tableau $T_d(v_{ a_1\dots a_l},
v_{ b_1\dots b_l})$ is nonmovable.
\end{prop}

\begin{pf}
Suppose the tableau $T_d(v_{ a_1\dots a_l},
v_{ b_1\dots b_l})$ in (\ref{eq:tskew})
is semistandard.
Then $ a_i< b_{i+d}$ holds for all $i=1,\dots,
l-d$. Thus for the permutation
\begin{equation*}
\sigma=\left(
\begin{array}{cccccc}
1&\dots&l-d&l-d+1&\dots&l\\
d+1&\dots&l&1&\dots&d
\end{array}
 \right),
\end{equation*}
$\displaystyle\sum_{i=1}^l H_1(v_{ a_i},
v_{ b_{\sigma(i)}})
=\displaystyle\sum_{i=l-d+1}^l H_1(v_{ a_i},
v_{ b_{\sigma(i)}})
\leq d$ holds. 
Therefore $H_l(v_{ a_1\dots a_l},v_{ b_1\dots b_l})\leq d$.

Conversely, suppose
$H_l(v_{ a_1\dots a_l},v_{ b_1\dots b_l})\leq d$. Then there is a
permutation $\sigma$ and a subset
$J$ of $\{ 1,2,\dots,l\}$ such that $\# J \geq l-d$
and $ a_i< b_{\sigma(i)}$ for each $i\in J$.
Now if $T_d
(v_{ a_1\dots a_l},v_{ b_1\dots b_l})$
is {\it not\/} semistandard, there exists
a number $j_0\leq l-d$ such that
$ a_{j_0}\geq  b_{j_0+d}$.
It follows that
\begin{equation*}
 b_1\leq  b_2\leq \cdots \leq
 b_{j_0+d}
\leq  a_{j_0}\leq
 a_{j_0+1}\leq
\cdots
\leq  a_l.
\end{equation*}
This implies $\#(J\cap
\{j_0,j_0+1,\dots,l\})\leq l-j_0-d$.
Thus we have $\#J\leq l-d-1$, which 
contradicts the assumption.
It is clear that (ii) is  equivalent to (i).
\end{pf}

\section{Spectral decomposition}\label{decom}

\subsection{Spectrum}

For $\vec{h}=(h_i)\in \{0,1,\dots,l\}^{\bold{N}}$
let $\tilde\kappa(\vec{h})$ be an
infinite skew diagram,
\begin{equation*}
\setlength{\unitlength}{0.8pt}
\begin{picture}(182,45)(-110,-40)
\put(-60,-40){\line(1,0){70}}
\put(-60,-30){\line(1,0){100}}
\put(-30,-20){\line(1,0){100}}
\put(0,-10){\line(1,0){70}}
\put(-60,-30){\line(0,-1){10}}
\put(-50,-30){\line(0,-1){10}}
\put(-40,-30){\line(0,-1){10}}
\put(-30,-20){\line(0,-1){10}}
\put(-20,-20){\line(0,-1){10}}
\put(-10,-20){\line(0,-1){10}}
\put(0,-10){\line(0,-1){10}}
\put(0,-30){\line(0,-1){10}}
\put(10,-10){\line(0,-1){10}}
\put(10,-30){\line(0,-1){10}}
\put(20,-10){\line(0,-1){10}}
\put(30,-20){\line(0,-1){10}}
\put(40,-20){\line(0,-1){10}}
\put(60,-10){\line(0,-1){10}}
\put(70,-10){\line(0,-1){10}}
\put(35,-5){\circle*{2}}
\put(40,0){\circle*{2}}
\put(45,5){\circle*{2}}
\put(-25,-35){\circle*{2}}
\put(-20,-35){\circle*{2}}
\put(-15,-35){\circle*{2}}
\put(5,-25){\circle*{2}}
\put(10,-25){\circle*{2}}
\put(15,-25){\circle*{2}}
\put(35,-15){\circle*{2}}
\put(40,-15){\circle*{2}}
\put(45,-15){\circle*{2}}
\put(-60,-29){$\overbrace{\hskip22pt}$}
\put(-50,-17){\small $h_1$}
\put(-30,-19){$\overbrace{\hskip22pt}$}
\put(-20,-7){\small $h_2$}
\put(-110,-20){$\tilde\kappa(\vec{h})=$}
\put(72,-20){,}
\end{picture}
\end{equation*}
where each row of $\tilde\kappa(\vec{h})$
has $l$ boxes.
For $\vec{h}\in \{0,1,\dots,l\}^{\bold{N}}$
and $\vec{s}=(s_i)\in B_l^{\bold{N}}$
with $s_i=(v_{ a_{i1}\dots
 a_{il}})$, we associate a tableaux 
$T_{\vec{h}}(\vec{s})$ of
shape $\tilde\kappa(\vec{h})$,
\begin{equation*}
\setlength{\unitlength}{1.5pt}
\begin{picture}(162,45)(-60,-30)
\put(-30,-30){\line(1,0){70}}
\put(-30,-20){\line(1,0){100}}
\put(0,-10){\line(1,0){100}}
\put(30,0){\line(1,0){70}}
\put(-30,-20){\line(0,-1){10}}
\put(-20,-20){\line(0,-1){10}}
\put(-10,-20){\line(0,-1){10}}
\put(0,-10){\line(0,-1){10}}
\put(10,-10){\line(0,-1){10}}
\put(20,-10){\line(0,-1){10}}
\put(30,0){\line(0,-1){10}}
\put(30,-20){\line(0,-1){10}}
\put(40,0){\line(0,-1){10}}
\put(40,-20){\line(0,-1){10}}
\put(50,0){\line(0,-1){10}}
\put(60,-10){\line(0,-1){10}}
\put(70,-10){\line(0,-1){10}}
\put(90,0){\line(0,-1){10}}
\put(100,0){\line(0,-1){10}}
\put(30.5,-10){\makebox(10,10){\small$ a_{31}$}}
\put(40.5,-10){\makebox(10,10){\small$ a_{32}$}}
\put(90.5,-10){\makebox(10,10){\small$ a_{3l}$}}
\put(0.5,-20){\makebox(10,10){\small$ a_{21}$}}
\put(10.5,-20){\makebox(10,10){\small$ a_{22}$}}
\put(60.5,-20){\makebox(10,10){\small$ a_{2l}$}}
\put(-29.5,-30){\makebox(10,10){\small$ a_{11}$}}
\put(-19.5,-30){\makebox(10,10){\small$ a_{12}$}}
\put(30.5,-30){\makebox(10,10){\small$ a_{1l}$}}
\put(65,5){\circle*{2}}
\put(70,10){\circle*{2}}
\put(75,15){\circle*{2}}
\put(5,-25){\circle*{2}}
\put(10,-25){\circle*{2}}
\put(15,-25){\circle*{2}}
\put(35,-15){\circle*{2}}
\put(40,-15){\circle*{2}}
\put(45,-15){\circle*{2}}
\put(65,-5){\circle*{2}}
\put(70,-5){\circle*{2}}
\put(75,-5){\circle*{2}}
\put(-30,-19){$\overbrace{\hskip40pt}$}
\put(-18,-13){\small $h_1$}
\put(0,-9){$\overbrace{\hskip40pt}$}
\put(12,-3){\small $h_2$}
\put(-60,-10){$T_{\vec{h}}(\vec{s})=$}
\put(102,-10){.}
\end{picture}
\end{equation*}

We introduce the {\it local energy map} $\rho:B_l^{\bold{N}}
\rightarrow \{ 0,1,\dots,l\}^{\bold{N}}$ as
\begin{equation}\label{eq:maph}
\rho:\vec{s}=(s_i)\mapsto \vec{h}=(h_i),\quad
h_i=H_l(s_{i+1},s_{i}).
\end{equation}

\begin{prop}\label{prop:diag}
(i) For any $\vec{s}\in B_l^{\bold{N}}$,
the tableau $T_{\rho(\vec{s})}(\vec{s})$
 is
nonmovable.\hfill\break
(ii) The map
\begin{equation*}
\tilde\varphi:\vec{s}\mapsto
T_{\rho(\vec{s})}(\vec{s}).
\end{equation*}
gives a one-to-one
correspondence between the elements
of $\rho^{-1}(\vec{h})$ and the
nonmovable tableaux of shape $\tilde\kappa(\vec{h})$.
\end{prop}
\begin{pf}
(i). Let $\vec{h}=\rho(\vec{s})$.
By definition, 
$T_{h_i}(v_{ a_{i+11}\dots a_{i+1l}},
v_{ a_{i1}\dots a_{il}})$
in (\ref{eq:tskew})
is the $i$th and $i+1$th rows (from the bottom) of
$T_{\vec{h}}(\vec{s})$.
Then it follows from Proposition \ref{prop:hbytableau}
 that $T_{\vec{h}}(\vec{s})$ is
nonmovable.
(ii). It is clear that $\tilde\varphi$
is injective. To see 
$\tilde\varphi$ maps
$\rho^{-1}(\vec{h})$ onto
the set of the nonmovable tableaux of shape
$\tilde\kappa(\vec{h})$,
suppose $T$ is a nonmovable
tableau of shape $\tilde\kappa(\vec{h})$.
Then there exists a unique $\vec{s}\in B_l^{\bold{N}}$
such that $T=T_{\vec{h}}(\vec{s})$.
{}From Proposition \ref{prop:hbytableau},
we have $\rho(\vec{s})=\vec{h}$.
\end{pf}

The {\it length\/} of a skew diagram
$\lambda/\mu$, denoted by  $l(\lambda/\mu)$, is defined as
$$
l(\lambda/\mu)=\max_i \{\lambda'_i-\mu'_i\},
$$
where $\lambda'$ is the conjugate of $\lambda$.

\begin{lem}\label{lem:imagecondition}
Let $\vec{h}=(h_i)\in \{0,1,\dots,l\}^{\bold{N}}$.
The following conditions are equivalent:
\begin{itemize}
\item[(i)]{$\vec{h}\in \rho(B_l^{\bold{N}})$.}
\item[(ii)]{There exists a nonmovable tableau of skew shape 
$\tilde\kappa(\vec{h})$.}
\item[(iii)]{The length of $\tilde\kappa(\vec{h})$
is at most $n$.}
\item[(iv)]{
$h_i+h_{i+1}+
\cdots +h_{i+n-1} \geq l
$
for any $i\geq 1$.}
\end{itemize}
\end{lem}
\begin{pf}
(i) $\Leftrightarrow$ (ii).
This is due to Proposition \ref{prop:diag}.

(ii) $\Leftrightarrow$ (iii).
Suppose the length of
$\tilde\kappa(\vec{h})$ is more than
 $n$. Then $\tilde\kappa(\vec{h})$ has
at least one column whose height is greater than
$n$. But, then, $\tilde\kappa(\vec{h})$ does not
admit any semistandard tableau
 of its shape (with $n$ numbers).
In particular, there is no nonmovable tableau
of shape $\tilde\kappa(\vec{h})$.
Conversely,
suppose the length of $\tilde\kappa(\vec{h})$
is at most $n$.
Then we can construct a nonmovable tableau of
shape $\tilde\kappa(\vec{h})$ by
filling the boxes of each column of
$\tilde\kappa(\vec{h})$
{}from the top to the bottom
 by the numbers
$1,2,3,\dots$ in this order.
It is clear that
the resulting tableau is semistandard.
Furthermore, the number in the top box of 
each column is always $1$. Such a tableau
is always nonmovable.

(iii) $\Leftrightarrow$ (iv).
This is clear from the shape of $\tilde\kappa(\vec{h})$.
\end{pf}

We define ${\mathrm{Sp}}_K=
\rho({\cal S}_K)$,
and call it
the {\it spectrum of\/} ${\cal S}_K$.
By definition $\mathrm{Sp}_K$ is a subset
of $\rho(B_l^{\bold{N}})$.
We set ${\mathrm{Sp}}=\bigsqcup_{K\in {\cal K}_l}
{\mathrm{Sp}}_K$.

\begin{lem}\label{lem:ground}
Let $K\in {\cal K}_l$.
\hfill\break
(i) It holds that
$\rho(\vec{s}^{(K)})=K$. In particular,
$K$ belongs to $\mathrm{Sp}_K$.
\hfill\break
(ii) For any $\vec{s}\in B_l^{\bold{N}}$,
$\vec{s}\approx \vec{s}^{(K)}$ if
and only if $\rho(\vec{s})\approx K$. 
\end{lem}
\begin{pf}
(i)
The tableau 
$T_{k_i}(s^{(K)}_{i+1}, s^{(K)}_i)$ is nonmovable as shown below:
\begin{equation*}
\setlength{\unitlength}{1pt}
\begin{picture}(192,35)(0,0)
\put(0,0){\line(1,0){160}}
\put(0,10){\line(1,0){190}}
\put(30,20){\line(1,0){160}}
\put(0,10){\line(0,-1){10}}
\put(30,20){\line(0,-1){20}}
\put(70,20){\line(0,-1){20}}
\put(130,20){\line(0,-1){20}}
\put(160,20){\line(0,-1){20}}
\put(190,20){\line(0,-1){10}}
\put(0,0){\makebox(30,10){\small $1\cdots 1$}}
\put(30,0){\makebox(40,10){\small $2\ \cdots\  2$}}
\put(130,0){\makebox(30,10){\small $n\cdots n$}}
\put(30,10){\makebox(40,10){\small $1\ \cdots\  1$}}
\put(160,10){\makebox(30,10){\small $n \cdots  n$}}
\put(0,11){\small $\overbrace{\hskip28pt}$}
\put(10,19){\small $k_{i}$}
\put(30,21){\small $\overbrace{\hskip38pt}$}
\put(48,29){\small $k_{i+1}$}
\put(192,0){.}
\end{picture}
\end{equation*}
Thus $H(s^{(K)}_{i+1}, s^{(K)}_i) = k_i$.
(ii)
If $\vec{s}
\approx \vec{s}^{(K)}$, then 
we have $\rho(\vec{s})\approx
\rho(\vec{s}^{(K)})=K$ by (i).
Conversely, let $\vec{h}=\rho(\vec{s})$,
and suppose $\vec{h}\approx K$.
Then except for finitely many
columns the heights of the columns of
$\tilde\kappa(\vec{h})$ are $n$.
The content
in any height-$n$ column
of a semistandard tableau 
is uniquely determined 
as $1,2,\dots,n$ from the top
to the bottom.
It follows that $\vec{s}\approx
\vec{s}^{(K)}$.
\end{pf}

{}From Lemmas \ref{lem:imagecondition}
and \ref{lem:ground} it follows that

\begin{thm}\label{prop:speccondition}
 An element $\vec{h}=(h_i)\in \{0,1,
\dots,l \}^{\bold{N}}$
belongs to ${\mathrm{Sp}}_K$ if and only if it satisfies
the conditions,
\begin{align}\label{eq:speccond}
\mbox{(i)} &\ \mbox{$h_i+h_{i+1}+
\cdots +h_{i+n-1} \geq l$ for any $i>0$.}
\tag{\theequation a}
\\
\mbox{(ii)} &\ %
\mbox{$\vec{h}\approx K$.}
\tag{\theequation b}
\end{align}
\addtocounter{equation}{1}
\end{thm}

\subsection{Spectral decomposition}

The local energy map $\rho$
 induces the {\it spectral
decomposition} of ${\cal S}_K$,
\begin{equation*}\label{eq:spindecomposition}
{\cal S}_K = \bigsqcup_{\vec{h}\in
{\mathrm{Sp}}_K} {\cal S}_{\vec{h}},
\quad {\cal S}_{\vec{h}}:=
\rho^{-1}(\vec{h}).
\end{equation*}

Let us introduce the character
of the degeneracy of the spectrum at $\vec{h}$,
\begin{equation}\label{eq:finitecharacter}
\chi_{\vec{h}}(x) =
\sum_{\vec{s}\in {\cal S}_{\vec{h}}}
e^{{\mathrm{wt}}(\vec{s})}.
\end{equation}
Due to Theorem \ref{thm:DJKMO},
the spectral decomposition of
${\cal S}_K$ induces the decomposition
of 
the character of ${\cal L}(\Lambda(K))$,
\begin{equation}\label{eq:fullcharacter}
{\mathrm{ch}}\, {\cal L}(\Lambda(K))(q,x)=
\sum_{\vec{h}\in {\mathrm{Sp}}_K}
q^{\sum_{i=1}^\infty i (h_i-k_i)} 
\chi_{\vec{h}}(x).
\end{equation}

\subsection{Character $\chi_{\vec{h}}$ and 
finite diagram}

Let us observe that for each element $\vec{h}=(h_i)
\in {\mathrm{Sp}}_K$, except for the case $\vec{h} = K$,
there is a unique positive integer $J$ such
that $\vec{h}$
 is
written in the following form
\begin{equation}\label{eq:finitepart}
\vec{h}=(h_1,\dots,h_J,
({k_{J+1}},\dots,
{k_{J+n}})^\infty),
\quad h_J > k_J.
\end{equation}
The subsequence $(h_1,\dots,h_{J})$
of $\vec{h}$ will be called the {\it finite part\/}
of $\vec{h}$, and denoted by
$\vec{h}_{\mathrm{fin}}$.
In the case
$\vec{h}=K$, its finite part is
defined to be
the empty sequence $\emptyset$.
For $\vec{h}\in {\mathrm{Sp}}_K$ with
$\vec{h}_{\mathrm{fin}}=(h_1,\dots,h_{J})$,
let $\kappa(\vec{h})=\kappa(h_1,\dots,h_J)$
be a finite skew subdiagram
of $\tilde\kappa(\vec{h})$
 as in Fig.\ \ref{fig:finite}.
The following properties and proposition are easily seen from
Fig.\ \ref{fig:finite}.
\begin{itemize}
\item[(i)]{The height of every column of
the compliment of $\kappa(\vec{h})$ in
$\tilde\kappa(\vec{h})$ is $n$.}
\item[(ii)]{The height of the rightmost column of
 $\kappa(\vec{h})$
is  at most $n-1$.}
\end{itemize}
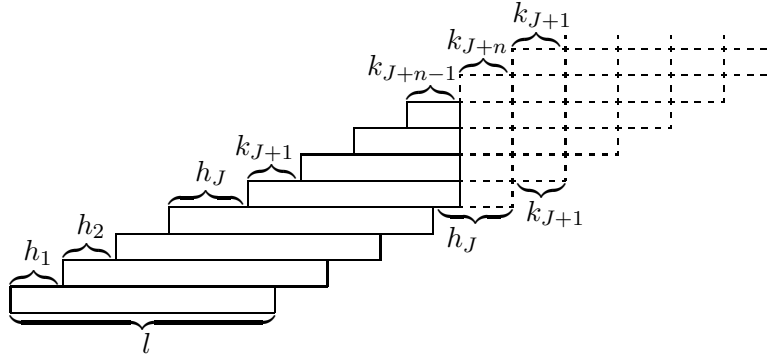
\begin{figure}[bt]
\begin{center}
\setlength{\unitlength}{1pt}
\begin{picture}(290,130)(-170,-115)
\put(-170,-100){\line(1,0){100}}
\put(-170,-90){\line(1,0){120}}
\put(-150,-80){\line(1,0){120}}
\put(-130,-70){\line(1,0){120}}
\put(-110,-60){\line(1,0){110}}
\put(-80,-50){\line(1,0){80}}
\put(-60,-40){\line(1,0){60}}
\put(-40,-30){\line(1,0){40}}
\put(-20,-20){\line(1,0){20}}
\put(0,-20){\line(0,-1){40}}
\put(-10,-60){\line(0,-1){10}}
\put(-20,-20){\line(0,-1){10}}
\put(-30,-70){\line(0,-1){10}}
\put(-40,-30){\line(0,-1){10}}
\put(-50,-80){\line(0,-1){10}}
\put(-60,-40){\line(0,-1){10}}
\put(-70,-90){\line(0,-1){10}}
\put(-80,-50){\line(0,-1){10}}
\put(-110,-60){\line(0,-1){10}}
\put(-130,-70){\line(0,-1){10}}
\put(-150,-80){\line(0,-1){10}}
\put(-170,-90){\line(0,-1){10}}
\multiput(0,-60)(5,0){4}{\line(1,0){1}}
\multiput(4,-60)(5,0){4}{\line(1,0){1}}
\multiput(0,-50)(5,0){8}{\line(1,0){1}}
\multiput(4,-50)(5,0){8}{\line(1,0){1}}
\multiput(0,-40)(5,0){12}{\line(1,0){1}}
\multiput(4,-40)(5,0){12}{\line(1,0){1}}
\multiput(0,-30)(5,0){16}{\line(1,0){1}}
\multiput(4,-30)(5,0){16}{\line(1,0){1}}
\multiput(0,-20)(5,0){20}{\line(1,0){1}}
\multiput(4,-20)(5,0){20}{\line(1,0){1}}
\multiput(0,-10)(5,0){24}{\line(1,0){1}}
\multiput(4,-10)(5,0){24}{\line(1,0){1}}
\multiput(20,0)(5,0){20}{\line(1,0){1}}
\multiput(24,0)(5,0){20}{\line(1,0){1}}
\multiput(0,-10)(0,-5){10}{\line(0,-1){1}}
\multiput(0,-14)(0,-5){10}{\line(0,-1){1}}
\multiput(20,0)(0,-5){12}{\line(0,-1){1}}
\multiput(20,-4)(0,-5){12}{\line(0,-1){1}}
\multiput(40,5)(0,-5){11}{\line(0,-1){1}}
\multiput(40,1)(0,-5){11}{\line(0,-1){1}}
\multiput(60,5)(0,-5){9}{\line(0,-1){1}}
\multiput(60,1)(0,-5){9}{\line(0,-1){1}}
\multiput(80,5)(0,-5){7}{\line(0,-1){1}}
\multiput(80,1)(0,-5){7}{\line(0,-1){1}}
\multiput(100,5)(0,-5){5}{\line(0,-1){1}}
\multiput(100,1)(0,-5){5}{\line(0,-1){1}}
\multiput(120,5)(0,-5){3}{\line(0,-1){1}}
\multiput(120,1)(0,-5){3}{\line(0,-1){1}}
\put(-170,-90){\small$\overbrace{\hskip0pt}$}
\put(-165,-80){\small$ h_1$}
\put(-150,-80){\small$\overbrace{\hskip0pt}$}
\put(-145,-70){\small$ h_2$}
\put(-110,-60){\small$\overbrace{\hskip28pt}$}
\put(-100,-50){\small$ h_J$}
\put(-80,-50){\small$\overbrace{\hskip0pt}$}
\put(-85,-40){\small$ k_{J+1}$}
\put(-20,-20){\small$\overbrace{\hskip0pt}$}
\put(-35,-10){\small$ k_{J+n-1}$}
\put(0,-10){\small$\overbrace{\hskip0pt}$}
\put(-5,0){\small$ k_{J+n}$}
\put(20,0){\small$\overbrace{\hskip0pt}$}
\put(20,10){\small$ k_{J+1}$}
\put(-170,-100){\small$\underbrace{\hskip100pt}$}
\put(-120,-115){\small$ l$}
\put(-8,-60){\small$\underbrace{\hskip28pt}$}
\put(-5,-75){\small$ h_{J}$}
\put(22,-50){\small$\underbrace{\hskip0pt}$}
\put(25,-65){\small$ k_{J+1}$}
\end{picture}
\end{center}
\caption{A finite diagram $\kappa(\vec{h})
=\kappa(h_1,\dots,h_J)$.}
\label{fig:finite}
\end{figure}

\begin{prop}\label{prop:bijection}
There is a one-to-one correspondence
between the nonmovable tableaux of
shape $\tilde\kappa(\vec{h})$ and those 
of shape $\kappa(\vec{h})$. The correspondence
is given by the restriction of a
nonmovable tableau of shape $\tilde\kappa
(\vec{h})$ on $\kappa(\vec{h})$.
\end{prop}

Combining the bijection in
Proposition \ref{prop:bijection} with the bijection
$\tilde\varphi$ in Proposition \ref{prop:diag},
we obtain a bijection
\begin{equation*}
\varphi: 
{\cal S}_{\vec{h}}
\mathop\rightarrow^\sim
{\mathrm{NMT}}(\kappa(\vec{h})).
\end{equation*}

Let us introduce the $sl_n$-weight of a tableau $T$ by
${\mathrm{wt}}(T) = \sum_{a=1}^n m_a \epsilon_a$,
where $(m_1, \ldots, m_n)$ is the content of $T$, i.e., 
$m_a$ is the number counting how many $a$'s are in $T$.

\begin{prop}\label{prop:weightpreserve}
The bijection $\varphi:
{\cal S}_{\vec{h}}
\rightarrow 
{\mathrm{NMT}}(\kappa(\vec{h}))$
 is weight-preserving,
i.e., for any $\vec{s}\in {\cal S}_{\vec{h}}$,
${\mathrm{wt}}(\vec{s})=
{\mathrm{wt}}(\varphi(\vec{s}))$ holds.
\end{prop}
\begin{pf}
Let $\vec{h}_{\mathrm{fin}}=(h_1,\dots,h_{J})$.
{}From Fig.\ \ref{fig:finite}, it is clear that
for any $\vec{s}=(s_i)\in \rho^{-1}(\vec{h})$,
 $s_i$ coincides with $s_i^{(K)}$ for any
$i\geq J+n$.
Thus, by definition,
\begin{equation*}
{\mathrm{wt}}(\vec{s})=
\sum_{i=1}^{n-1} k_i \overline\Lambda_{i}
+ \sum_{i=1}^{J+n-1} {\mathrm{wt}}(s_i)
- \sum_{i=1}^{J+n-1} {\mathrm{wt}}(s_i^{(K)}).
\end{equation*}
In the meanwhile
\begin{align*}
\sum_{i=1}^{J+n-1} {\mathrm{wt}}(s_i)
&={\mathrm{wt}}(\varphi(\vec{s}))
-\sum_{i=1}^{n-1} k_{J+i-1}
\overline\Lambda_{i},\\
\sum_{i=1}^{J+n-1} {\mathrm{wt}}(s_i^{(K)})
&=
\sum_{i=1}^{n-1} k_{i}
\overline\Lambda_{i}
-\sum_{i=1}^{n-1} k_{J+i-1}
\overline\Lambda_{i}.
\end{align*}
Therefore we have ${\mathrm{wt}}(\vec{s})=
{\mathrm{wt}}(\varphi(\vec{s}))$.
\end{pf}

We define the function
\begin{equation}
t_{\lambda/\mu}(x) =
\sum_{T\in{\mathrm{NMT}}(\lambda/\mu)}
e^{{\mathrm{wt}}(T)}.
\end{equation}
Let ${\mathrm{Sp}}_{K,\mathrm{fin}}$ denote the
set of all the finite parts of $\vec{h}$'s in 
${\mathrm{Sp}}_K$, i.e.,
\begin{gather*}
{\mathrm{Sp}}_{K,\mathrm{fin}}=
\{(h_1,\dots,h_{J}) \mid
J\geq0,\ h_i\in \{ 0,1,\dots,l\},\ 
\mbox{the condition ($\ast$) when $J\geq 1$}
\},
\end{gather*}
where $(h_1,\dots,h_J)=\emptyset$ for
$J=0$, and the condition ($\ast$) is
\begin{gather*}
(\ast)
\begin{cases}
& \sum_{j=0}^{n-1}h_{i+j}\geq l
\quad\mbox{for any $1\le i \le J-1$,}
\quad \mbox{($h_{J+i}:={k_{J+i}}$,\ 
$1\le i \le n-1$)},\cr
& h_J>{k_{J}}.
\end{cases}
\end{gather*}
{}From Proposition 
\ref{prop:weightpreserve} and
(\ref{eq:fullcharacter}), we have
\begin{thm}\label{thm:maintheorem}
(i) The character $\chi_{\vec{h}}$ of
${\cal S}_{\vec{h}}$ is equal to the function
$t_{\kappa(\vec{h})}.$
\hfill\break
(ii)
The character of the level $l$ integrable module
${\cal L}(\Lambda(K))$ of $\widehat{sl}_n$
decomposes as
\begin{equation*}\label{eq:characterdecomposition}
\begin{split}
{\mathrm{ch}}\, {\cal L}(\Lambda(K))(q,x)
&=
\sum_{\vec{h}\in
{\mathrm{Sp}}_{K}}
q^{\sum_{i=1}^\infty i (h_i - k_{i})}
t_{\kappa(\vec{h})}(x)\\
&=
\sum_{(h_1,\dots,h_{J})\in
{\mathrm{Sp}}_{K,\mathrm{fin}}}
q^{\sum_{i=1}^Ji (h_i - k_{i})}
t_{\kappa(h_1,\dots,h_J)}(x)
\end{split}
\end{equation*}
\end{thm}

\subsection{Formula of 
$t_\kappa(h_1,\dots,h_J)$ by skew Schur functions}

Let 
$s_{\lambda/\mu} = s_{\lambda/\mu}(x)$ 
be the {\it skew Schur function associated to} $\lambda/\mu$ \cite{Ma}.
Let $(h_1,\dots,h_J)\in {\mathrm{Sp}}_{K,\mathrm{fin}}$
be fixed.
We set
$I_0:=\{ i\in \{1,\dots,J\}|h_i\neq 0\}$ and
$\kappa_0:=\kappa(h_1,\dots,h_J)$.
For  
$i_1,\dots,i_p\in I_0$ such that
$i_1<i_2<\dots < i_p$, 
consider a sequence $(h'_1,\dots,h'_J)$, where
\begin{equation*}
h'_i=
\begin{cases}
h_i-1 & \mbox{if $i=i_1,\dots,i_p$,}\\
h_i & \mbox{otherwise.}
\end{cases}
\end{equation*}
Even though 
$(h'_1,\dots,h'_J)$
does not necessarily
belong to ${\mathrm{Sp}}_{K,\mathrm{fin}}$,
we can still define a skew diagram $\kappa(h'_1,
\dots,h'_J)$ by Fig.\ \ref{fig:finite}.
Let 
$\kappa_{i_1,\dots,i_p}:=
\kappa(h'_1,\dots,h'_J)$.

\begin{prop}\label{prop:characterformula}
Let $(h_1,\dots,h_J)\in
{\mathrm{Sp}}_{K,\mathrm{fin}}$. Then
\begin{equation}\label{eq:altsum}
t_{\kappa(h_1,\dots,h_J)}
=
s_{\kappa_0}
+\sum_{p=1}^{|I_0|}
(-1)^p
\sum_{i_1,\dots,i_p\in I_0\atop i_1<\cdots<i_p}
s_{\kappa_{i_1,\dots,i_p}}.
\end{equation}
where in the right hand side we impose the relation
$x_1 \cdots x_n= 1$.
\end{prop}

\begin{exm}
\begin{equation}
\setlength{\unitlength}{1pt}
\begin{picture}(173,17)
\put(0,0){\line(1,0){10}}
\put(0,5){\line(1,0){15}}
\put(5,10){\line(1,0){15}}
\put(10,15){\line(1,0){10}}
\put(0,0){\line(0,1){5}}
\put(5,0){\line(0,1){10}}
\put(10,0){\line(0,1){15}}
\put(15,5){\line(0,1){10}}
\put(20,10){\line(0,1){5}}
\put(-1,11){$t$}
\put(24,8){$ = $}
\put(39,11){$s$}
\put(40,0){\line(1,0){10}}
\put(40,5){\line(1,0){15}}
\put(45,10){\line(1,0){15}}
\put(50,15){\line(1,0){10}}
\put(40,0){\line(0,1){5}}
\put(45,0){\line(0,1){10}}
\put(50,0){\line(0,1){15}}
\put(55,5){\line(0,1){10}}
\put(60,10){\line(0,1){5}}
\put(80,0){\line(1,0){10}}
\put(80,5){\line(1,0){15}}
\put(85,10){\line(1,0){10}}
\put(85,15){\line(1,0){10}}
\put(80,0){\line(0,1){5}}
\put(85,0){\line(0,1){15}}
\put(90,0){\line(0,1){15}}
\put(95,5){\line(0,1){10}}
\put(64,8){$ - $}
\put(77,11){$s$}
\put(120,0){\line(1,0){10}}
\put(120,5){\line(1,0){10}}
\put(120,10){\line(1,0){15}}
\put(125,15){\line(1,0){10}}
\put(120,0){\line(0,1){10}}
\put(125,0){\line(0,1){15}}
\put(130,0){\line(0,1){15}}
\put(135,10){\line(0,1){5}}
\put(102,8){$ - $}
\put(114,11){$s$}
\put(160,0){\line(1,0){10}}
\put(160,5){\line(1,0){10}}
\put(160,10){\line(1,0){10}}
\put(160,15){\line(1,0){10}}
\put(160,0){\line(0,1){15}}
\put(165,0){\line(0,1){15}}
\put(170,0){\line(0,1){15}}
\put(152,11){$s$}
\put(139,8){$+$}
\put(172,8){$.$}
\end{picture}
\nonumber
\end{equation}
\end{exm}

\begin{pf}
For a given skew shape $\kappa$, let us denote by SST$(\kappa)$ the set
of all semistandard tableaux of the shape $\kappa$. Consider the injection 
\begin{equation*}
\eta_{i_1,\dots,i_p}:
{\mathrm{SST}}(\kappa_{i_1,\dots,i_p})
\hookrightarrow 
{\mathrm{SST}}(\kappa_{0})
\end{equation*}
preserving the contents for each row.
It is clear that the map
$\eta_{i_1,\dots,i_p}$ is weight-preserving.
Let $A_{i_1,\dots,i_p}$ denote the image
of ${\mathrm{SST}}(\kappa_{i_1,\dots,i_p})$
in ${\mathrm{SST}}(\kappa_{0})$ under the map.
Then we have 
\begin{itemize}
\item[(i)]{${\mathrm{NMT}}(\kappa_0)=
{\mathrm{SST}}(\kappa_0)-\bigcup_{i\in I_0}A_i$.}
\item[(ii)]{$A_{i_1,\dots,i_p}=A_{i_1}\cap
A_{i_2}\cap\dots \cap A_{i_p}$.}
\end{itemize}
{}From (i),(ii) and the inclusion-exclusion principle
we obtain (\ref{eq:altsum}).
\end{pf}


In the case of either
$l=1$ or $n=2$, the right hand side of
the formula (\ref{eq:altsum}) is written
in a simpler form,
 and reproduces the known results.

\begin{exm}
{\it The case $l=1$, $n$: general (cf.\ \cite{KKN}).}
The skew diagram $\kappa_0=
\kappa(h_1,\dots,h_J)$ has $I_0+1$ columns,
and the width of every row of $\kappa_0$ is $1$.
Let $1\leq m_i$ $(i=1,\dots,I_0+1)$ be the height of the
$i$th column (from the left). Let $E_m$ be the
$m$th elementary symmetric polynomial of $x_1,\dots,x_n$
for $m=0,1,\dots,n$, and $E_m=0$ for $m<0$ and $m>n$.
Then
\begin{equation*}
s_{\kappa_0}=
\prod_{i=1}^{I_0+1} E_{m_i},
\qquad
s_{\kappa_{i_1}}=
\left(\prod_{i=1}^{j_1-1} E_{m_i}\right)
E_{m_{j_1}+m_{j_1+1}}
\left(\prod_{i=j_1+2}^{I_0+1} E_{m_i}\right),
\qquad\dots,
\end{equation*}
where $j_1$ is the number determined by
$i_1=m_1+m_2+\cdots + m_{j_1}$. Thus the rhs of
(\ref{eq:altsum}) is written as
\begin{equation*}
\det_{1\leq i,j\leq I_0+1} E_{(\sum_{t=1}^j m_t)-
(\sum_{t=1}^{i-1}m_t)}.
\end{equation*}
Let $[m_1,\dots,m_{I_0+1}]$ be the {\it border strip}
with $I_0+1$ columns such that the $i$th
column (from the left) has height $m_i$.
(A border strip is a connected skew diagram
containing no $2\times 2$ block of boxes.)
Due to the Jacobi-Trudi formula, the above determinant
is equal to the skew Schur 
function $s_{[m_1,\dots,m_{I_0+1}]}$.
\end{exm}

\begin{exm}
{\it The case $n=2$, $l$: general (cf.\ \cite{arakawa}).}
The skew diagram $\kappa_0=
\kappa(h_1,\dots,h_J)$ has $J+1$ rows.
Let $m_1=h_1$ and $m_i=h_i+h_{i-1}-l$ for
$i=2,\dots,J+1$.
In other words, $m_i$ is the number of the boxes
in
the $i$th row of $\kappa_0$
without any vertical adjacent boxes.
For the simplicity, let us first consider the
situation $m_i>0$ for all $i$.
Let $H_m$ be the $m$th completely symmetric
polynomial of $x_1$ and $x_2$.
Then 
\begin{equation*}
s_{\kappa_0}=
\prod_{i=1}^{J+1} H_{m_i},
\qquad
s_{\kappa_{i_1}}=
\left(\prod_{i=1}^{i_1-1} H_{m_i}\right)
H_{m_{i_1}-1} H_{m_{i_1+1}-1}
\left(\prod_{i=i_1+2}^{J+1} H_{m_i}\right),
\qquad\dots 
\end{equation*}
By the formulae
\begin{align*}
H_{m_1+m_2}&=H_{m_1}H_{m_2}-H_{m_1-1}H_{m_2-1},\cr
H_{m_1+m_2+m_3}&=H_{m_1}H_{m_2}H_{m_3}
- H_{m_1-1}H_{m_2-1}H_{m_3}\cr
&\qquad - H_{m_1}H_{m_2-1}H_{m_3-1}
+ H_{m_1-1}H_{m_2-2}H_{m_3-1},
\qquad
\mbox{etc.},
\end{align*}
the rhs of (\ref{eq:altsum}) reduces to
$H_{m_1+\dots + m_{J+1}}$.
In general the sequence $(m_i)$ is decomposed
into $r\geq 1$ components of successive
nonzero elements as
\begin{equation*}
(m_i)=(
m_{11},\dots,m_{1t_1},
0,\dots,0,
m_{21},\dots,m_{2t_2},
0,\dots,0,
m_{r1},\dots,m_{rt_r}),
\qquad m_{ij}\neq 0.
\end{equation*}
Let $M_i=\sum_{j=1}^{t_i} m_{ij}$.
Then the rhs of (\ref{eq:altsum}) is
equal to $\prod_{i=1}^r H_{M_i}$.
\end{exm}


\section{Spectral decomposition and exponents}\label{exp}

Let $\lambda$ be a partition of length $\le m$, $\mu$ be a composition, 
$l(\mu )\le m$, and $d\in {\mathbf{Z}}^{m-1}_{\ge 0}$ be a sequence of 
nonnegative integers. In this section we are going to describe a bijection
\begin{equation}
{\mathrm{SST}}_d(\lambda ,\mu )\leftrightarrow 
{\mathrm{LR}}_0({\mathrm{Sh}}_d(\mu 
),\lambda ) \label{5.1}
\end{equation}
between the set of all the semistandard Young tableaux $T$ of shape $\lambda$ 
and content $\mu$, with a given set of exponents $d(T)=d=(d_1,d_2,\ldots 
)$, and the set of all the nonmovable Littlewood-Richardson 
tableaux of (skew) shape ${\mathrm{Sh}}_d(\mu )$ and content $\lambda$. We 
start with reminding all necessary definitions.

\subsection{Exponents} (\cite{GZ,KR,KB}). Let us 
explain at first a combinatorial definition of the exponents. We start 
with definition of the descent set $D(T)$ of a given semistandard tableau 
$T\in \mathrm{SST}(\lambda ,\mu )$ of shape $\lambda$ and content $\mu$. Let 
indices $i$ and $i+1$ belong to the given tableau $T$. We say that a pair 
$i$ and $i+1$ forms a descent in tableau $T$, if $i+1$ lies strictly 
below than $i$ in the tableau $T$. We say that $i$ is a descent of 
multiplicity $\zeta_i$, if $\zeta_i$ is the maximal number such that 
there exist $\zeta_i$ pairs of descents $i$ and $i+1$ in the tableau $T$ 
with different ends $i+1$. We denote by $D(T)$ a set of all descents with 
multiplicities in the tableau $T$.

\begin{exm}\label{e5.1} Let us take
$$T:=\begin{array}{llllll} 1&1&2&3&5&5\\ 2&3&3&4\\ 3&4&5\\ 4&5 
\end{array}\in \mathrm{SST}((6,4,3,2), (2,2,4,3,4)),
$$
then $D(T)=\{ 1,2,2,3,3,3,4,4\}$ and the multiplicities of descents are 
$\zeta =(1,2,3,2)$.
\end{exm}

\begin{defn}\label{d5.1} We define the $i$th exponent $d_i(T)$ of a 
semistandard tableau $T\in \mathrm{SST}(\lambda ,\mu )$ as follows
$$d_i(T):=\mu_{i+1}-\zeta_i(T), \ \ 1\le i\le m-1.
$$
\end{defn}
In the Example~\ref{e5.1} we have $\mu =(2,2,4,3,4)$ and $d=(1,2,0,2)$.

There exists also a group-representation interpretation of the 
exponents (cf. \cite{GZ}):
$$d_i(T)=\max_k\{ E_{i,i+1}^k~|~T\rangle\ne 0\} ,
$$
where $E_{i,j}$, $1\le i\ne j\le m$, is a set of generators of the Lie 
algebra $\hbox{\germ gl}_m$, and $|T\rangle$ is the Gelfand-Tsetlin pattern 
(GT-pattern, for short) corresponding to the tableau $T$. Here we 
identify the GT-patterns with a basis in the highest weight $\lambda$ 
irreducible representation of $\hbox{\germ gl}_m$.

\subsection{Littlewood-Richardson tableaux.}\label{lrtab}

Let $A$ be a skew shape and $\mu$ be a composition, $|A|=|\mu |$, $l(\mu 
)=m$. We denote by ${\mathrm{Tab}}(A,\mu )$ 
the set of all the tableaux of shape $A$ 
and content $\mu$, in other words, the set of all fillings of the shape 
$A$ by the numbers from 1 to $m$ that have the content $\mu$.

\begin{exm}\label{e5.2} Let us take $\mu =(8,5,2)$ and 
$A=(6,6,6,6)/(4,3,2)$. Then
$$\begin{array}{llllll} &&&&1&1\\ &&&1&2&2\\ &&1&1&1&3\\ 1&1&2&2&2&3 
\end{array}\in{\mathrm{Tab}}(A,\mu ).
$$
\end{exm}

For each tableau $T$ we define a word $w(T)$ by reading the numbers of 
the tableau $T$ along the rows from the right to the left and from the 
top to the bottom. For the tableau from Example~\ref{e5.2}, we have
$$w(T)=112213111322211.
$$
\begin{defn}\label{little}
Let us say that a tableau $T$ is the Littlewood-Richardson 
tableau (LR-tableau, for short), if the corresponding word $w(T)$ is a 
lattice word and each row of $T$ is 
semistandard.
\end{defn}

For the reader's convenience, let us remind a definition of a lattice word 
(see, e.g. \cite{Ma}, p.143):

A word $w=a_1a_2\ldots a_N$ in the symbols $1,2,\ldots ,m$ is said to be 
a lattice word if for $1\le r\le N$ and $1\le i\le m-1$, the number of 
occurrences of the symbol $i$ in $a_1a_2\ldots a_r$ is not less than the 
number of occurrences of $i+1$.

In Example~\ref{e5.2}, the tableau $T$ is an LR-tableau.

\subsection{Bijection $\theta_{d}$} Let $\mu =(\mu_1,\ldots ,\mu_m )$ 
and $\nu =(\nu_1,\ldots ,\nu_{m-1})$ be compositions. We define a skew 
shape ${\rm Sh}_{\nu}(\mu )$ as
\begin{equation*}
\begin{picture}(204,90)(-62,-15)
\put(0,0){\line(1,0){60}}
\put(0,10){\line(1,0){80}}
\put(20,20){\line(1,0){60}}
\put(60,40){\line(1,0){60}}
\put(60,50){\line(1,0){80}}
\put(80,60){\line(1,0){60}}
\put(0,10){\line(0,-1){10}}
\put(60,10){\line(0,-1){10}}
\put(20,20){\line(0,-1){10}}
\put(80,20){\line(0,-1){10}}
\put(60,50){\line(0,-1){10}}
\put(120,50){\line(0,-1){10}}
\put(80,60){\line(0,-1){10}}
\put(140,60){\line(0,-1){10}}
\put(0,10){\small$\overbrace{\hskip0pt}$}
\put(-5,20){\small$ \nu_{m-1}$}
\put(60,50){\small$\overbrace{\hskip0pt}$}
\put(65,60){\small$ \nu_{1}$}
\put(80,60){\small$\overbrace{\hskip60pt}$}
\put(105,70){\small$ \mu_{1}$}
\put(0,0){\small$\underbrace{\hskip60pt}$}
\put(27,-15){\small$ \mu_{m}$}
\put(20,20){\small$\overbrace{\hskip60pt}$}
\put(42,30){\small$ \mu_{m-1}$}
\put(80,35){\circle*{2}}
\put(75,30){\circle*{2}}
\put(-62,30){${\mathrm{Sh}}_{\nu}(\mu )=$}
\put(142,30){.}
\end{picture}
\end{equation*}
In the sequel we will assume that $\nu_i+\mu_i\ge\mu_{i+1}$ holds
for any $i$. 
Now we are ready to define a map
\begin{equation}
\theta_{\nu}: \mathrm{SST}(\lambda ,\mu )\to
{\mathrm{Tab}}({\mathrm{Sh}}_{\nu}(\mu 
),\lambda ) \label{5.2}
\end{equation}
in the following way: Consider a given semistandard tableau 
$T\in \mathrm{SST}(\lambda ,\mu )$, 
let us fill the $k$th row of the shape ${\rm 
Sh}_{\nu}(\mu )$ by numbers $1,2,\ldots $, according to the indices of 
the rows in the tableau $T$ which contain $k$, starting from the top to 
the bottom.
\begin{exm}\label{e5.3} Let us take
$$T=\begin{array}{llllllll}1&1&2&3&3&3&4&4\\ 2&2&4&4&4\\ 3&4 \end{array},
$$
and $\nu =(1,1,2)$. Then
$$\theta_{\nu}(T):=\begin{array}{llllll} &&&&1&1\\ &&&1&2&2\\ &&1&1&1&3\\ 
1&1&2&2&2&3 \end{array}.
$$
\end{exm}

The following Lemma is well-known.
\begin{lem}\label{l5.1} Let $T$ be a semistandard tableau of shape 
$\lambda$, content $\mu$, and $\nu\in{\bf Z}^{m-1}_{\ge 0}$ be a 
composition. Then $\theta_{\nu}(T)$ is an LR-tableau, and $\theta_{\nu}$ 
defines a bijection
\begin{equation}
\theta_{\nu}: \mathrm{SST}(\lambda ,\mu )\leftrightarrow 
\mathrm{LR}({\mathrm{Sh}}_{\nu}(\mu ),\lambda ), \label{5.3}
\end{equation}
where ${\mathrm{LR}}(A, \eta )$ is the set of all the LR-tableaux 
(not necessarily semistandard!) of shape $A$ and content $\eta$.
\end{lem}
In other words, a tableau $T$ is semistandard if and only if for some 
composition $\nu\in{\mathbf{Z}}_{\ge 0}^{m-1}$, the tableau $\theta_{\nu}(T)$ 
is an LR-tableau.
It is easy to see that in general, a tableau $\theta_{\nu}(T)$, 
$T\in \mathrm{SST}(\lambda ,\mu )$, is not necessarily a semistandard one.
\begin{lem}\label{l5.2} Let $T$ be a semistandard tableau of shape 
$\lambda$ and content $\mu$, and $\nu\in {\mathbf{Z}}_{\ge 0}^{m-1}$ be a 
composition satisfying $\nu_i + \mu_i \ge \mu_{i+1}$. 
Then the tableau $\theta_{\nu}(T)$ is semistandard if and only if
\begin{equation}
\nu_i\ge d_i(T)
\qquad
\mbox{for any\/ $1\le i\le m-1$}. \label{5.4}
\end{equation}
\end{lem}
\begin{pf}
Suppose $\theta_\nu(T)$ is semistandard.
The overlapping length between the $i$th and the $i+1$th 
rows in $\theta_\nu(T)$ is $s =\mu_{i+1} - \nu_i$.
Let $k_1, \ldots, k_s$ ($m_1, \ldots, m_s$) be the 
entries of the first (last) $s$ boxes in the $i$th  ($i+1$th) row.
For each pair $(k_j,m_j)$ one can assign a descent pair
in $T$ with $i \ (i+1)$ at the $k_j$th ($m_j$th) row therein.
This means $\zeta_i \ge s = \mu_{i+1} - \nu_i$.
The converse is similar.
\end{pf}

It follows from Lemma~\ref{l5.2} that
\begin{cor}\label{c5.1} For a given semistandard tableau $T$ and a 
composition $\nu\in{\bf Z}_{\ge 0}^{m-1}$, the tableau $\theta_{\nu}(T)$ 
is nonmovable (see Definition~\ref{d1.1}) if and only if
$$\nu_i=d_i(T)\qquad
\mbox{for any\/ $1\le i\le m-1$}.
$$
\end{cor}

Before stating our main result of Section~\ref{exp}, let us introduce two 
additional notations. 

The first one is:~~
Given a partition $\lambda$, $l(\lambda )\le m$, and compositions $\mu$, 
$l(\mu )\le m$, and $d$, $l(d)\le m-1$, we denote by ${\mathrm{SST}}_d(\lambda 
,\mu )$ the set of all the semistandard tableaux $T$ of shape $\lambda$ and 
content $\mu$ such that $d_i(T)=d_i$, $1\le i\le m-1$. 
The second one is:~~
Given a skew shape $A$ and a partition $\lambda$, we denote by 
${\mathrm{LR}}_0(A,\lambda )$ the set of all the nonmovable 
LR-tableaux of (skew) shape 
$A$ and content $\lambda$.
See Definitions ~\ref{d1.1} and \ref{little}.

Combining together Lemmas~\ref{l5.1}, \ref{l5.2} and 
Corollary~\ref{c5.1}, we obtain the following result which is the main 
result of Section~\ref{exp}.
\begin{thm}\label{t5.1} Let $\lambda ,\mu$ and $d$ as above, then the map 
$\theta_d$ defines a bijection
\begin{equation}
\theta_d: {\mathrm{SST}}_d(\lambda ,\mu )
\leftrightarrow{\mathrm{LR}}_0({\mathrm{Sh}}_d(\mu ),\lambda ). \label{5.5}
\end{equation}
\end{thm}
\begin{cor}\label{c5.2} {\rm (``Spectral decomposition" of the Kostka 
numbers)}. Let $\lambda$ be a partition, and $\mu$ be a composition. Then
\begin{equation}
|{\rm SST}(\lambda ,\mu )|=\dim V_{\lambda}(\mu )=K_{\lambda ,\mu}=
\sum_{d=(d_1,\ldots ,d_{m-1})\in{\bf Z}_{\ge 0}^{m-1}}
|{\mathrm{LR}}_0({\mathrm{Sh}}_d(\mu ),\lambda )|. \label{5.6}
\end{equation}
\end{cor}

Let us 
say a few words about the set ${\mathrm{LR}}_0(\alpha /\beta ,\lambda )$.
Let $c_{\lambda\mu}^{\nu}$ be the Littlewood-Richardson coefficient, that is 
the multiplicity of the irreducible highest weight $\nu$ representation 
$V_{\nu}$ of the Lie algebra $\hbox{\germ gl}_m$ in the tensor product 
$V_{\lambda}\otimes V_{\mu}$.
\begin{prop}\label{p5.2} Let $\alpha ,\beta$ be partitions, 
$\beta\subset\alpha$. Then
$$|{\mathrm{LR}}_0(\alpha /\beta ,\lambda 
)|=\sum_{\alpha'/\beta'}(-1)^{{\rm sign}(\alpha /\beta ,\alpha'/\beta')}
c_{\beta'\lambda}^{\alpha'},
$$
where the sum extends over all the skew diagrams $\alpha'/\beta'$ such that 
$\alpha'_i-\beta'_i=\alpha_i-\beta_i$, $\beta'_i=\beta_i$, or $\beta_i-1$ 
for all $i$, and ${\rm sign}(\alpha /\beta ,\alpha'/\beta'):=|\{ 
i~|~\beta'_i=\beta_i-1\} |$.
\end{prop}
\begin{pf} There is a one-to-one correspondence between the 
semistandard LR-tableaux of shape $\alpha /\beta$ and content $\lambda$ 
which are movable to those of shape $\alpha'/\beta'$, and the 
semistandard LR-tableaux of shape $\alpha'/\beta'$ and content 
$\lambda$. It is well-known (see, e.g. \cite{Ma}, p.143) that the number 
of the last tableaux is equal to the Littlewood-Richardson coefficient 
$c_{\beta'\lambda}^{\alpha'}$. To finish the proof, let us apply the 
inclusion-exclusion principle.
\end{pf}

\begin{cor}\label{c5.3}  Let $A$ be a skew diagram. Let us denote by $\widehat 
A$ the diagram which is obtained from the diagram $A$ by the 180 
degree rotation. Then
\begin{equation}
|{\mathrm{LR}}_0(A,\nu )|=|{\mathrm{LR}}_0({\widehat A},\nu )|. \label{5.7}
\end{equation}
\end{cor}

\begin{pf} Assume that $A=\lambda /\mu$, where $\lambda 
=(\lambda_1,\ldots ,\lambda_m)$ and $\mu =(\mu_1,\ldots ,\mu_m)$. It is 
easy to see that ${\widehat A}={\widehat\mu}/{\widehat\lambda}$, where 
${\widehat\lambda}=(\lambda_1-\lambda_m,\lambda_1-\lambda_{m-1},\ldots 
,\lambda_1-\lambda_2,0)$ and 
${\widehat\mu}=(\lambda_1-\mu_m,\lambda_1-\mu_{m-1},\ldots ,\lambda_1-\mu_1)$. 
Now we are going to use the following well-known (see, e.g. \cite{KB}) 
symmetry property of the Littlewood-Richardson coefficients:
\begin{equation}
c_{\mu\nu}^{\lambda}=c_{{\widehat\lambda}\nu}^{\widehat\mu}. \label{5.8}
\end{equation}
{}From (\ref{5.8}) and Proposition~\ref{p5.2} we obtain the equality 
(\ref{5.7}).
\end{pf}

In Section~\ref{kostka} 
we are going to describe a $q$-analog of (\ref{5.6}) in 
a particular case when $\mu$ is a rectangular partition $(l^m)$.

\section{New combinatorial formula for Kostka--Foulkes
polynomials}\label{kostka}

In Section~\ref{kostka} we give a new combinatorial formula for the
Kostka-Foulkes 
polynomials $K_{\lambda ,\mu}(q)$ in the case when $\mu$ is a rectangular 
partition.
For the reader's convenience, we remind the basic definitions and 
results concerning the Kostka-Foulkes polynomials. The proofs and 
further details can be found in \cite{Ma}.

\subsection{Kostka-Foulkes polynomials.} Let $\lambda$ and $\nu$ 
be partitions such that $\nu_i\le\lambda_i$ for any $i$, and 
$l(\lambda/\nu )\le m$. 
Let $x=(x_1,\ldots ,x_m)$ be set of independent 
variables. The connection coefficients between the skew Schur functions 
$s_{\lambda /\nu}(x)$ and the monomial symmetric functions $m_{\mu 
}(x)$ are called the (skew) Kostka numbers:
\begin{equation}
s_{\lambda /\nu}(x)=\sum_{\mu}K_{\lambda /\nu ,\mu}m_{\mu}(x).
\label{6.2}
\end{equation}
It is well-known that the Kostka number $K_{\lambda /\nu ,\mu}$ 
is equal to the number of all the semistandard tableaux of shape 
$\lambda /\nu$ and content $\mu$.
\begin{defn}\label{d6.2} The Hall-Littlewood function $P_{\lambda}(x;q)$ 
corresponding to a partition $\lambda$ is defined by the following formula
\begin{equation}
P_{\lambda}(x;q)=c_{\lambda ,m}(q)\sum_{w\in S_m}w\left( x^{\lambda}
\prod_{i<j}\frac{x_i-qx_j}{x_i-x_j}\right) ,\label{6.3}
\end{equation}
where $c_{\lambda ,m}(q)$ is a normalization constant defined by 
\begin{eqnarray*}
&&c_{\lambda ,m}(q)=v_{m-l(\lambda )}(q)/v_{\lambda}(q), \\
&&v_{\lambda}(q)=\prod_{i\ge 1}\prod_{j=1}^{\lambda'_i-\lambda'_{i+1}}
\frac{1-q^j}{1-q}.
\end{eqnarray*}
\end{defn}
\begin{defn}\label{d6.3} The Kostka-Foulkes polynomials 
$K_{\lambda /\nu ,\mu}(q)$ are defined as the connection 
coefficients between the skew Schur and Hall-Littlewood functions:
\begin{equation}
s_{\lambda /\nu}(x)=\sum_{\mu}K_{\lambda /\nu ,\mu}(q)
P_{\mu}(x;q). \label{6.4}
\end{equation}
\end{defn}

A combinatorial description of the Kostka-Foulkes polynomials 
$K_{\lambda ,\mu}(q)$ has been found by A.~Lascoux and 
M.-P.~Sch\"utzenberger \cite{LS}, and goes as follows. Given two 
partitions $\lambda$ and $\mu$, it is possible to attach to each 
semistandard tableau $T\in{\rm SST}(\lambda ,\mu )$ a positive integer 
$c(T)$ (charge of tableau $T$) such that
\begin{equation}
K_{\lambda ,\mu}(q)=\sum_{T \in{\mathrm{SST}}(\lambda ,\mu )}q^{c(T)}.
\label{6.5}
\end{equation}
Below we give a definition of the charge of a tableau, 
according to  A.~Lascoux and M.-P.~Sch\"utzenberger, \cite{LS,Ma}. 
Let $\lambda$ and $\mu$ be partitions, and let $T\in{\rm 
SST}(\lambda ,\mu)$. Consider the word $w(T)$ which corresponds to the 
tableaux $T$, see \cite{Ma}, Chapter~I, \S 9, or our Section~\ref{lrtab}. The 
charge $c(T)$ of the tableau $T$ is defined as the charge of 
corresponding word $w(T)$. Now we define the charge of a word $w$. Recall 
that the weight $\mu$ of a word $w$ is a sequence $\mu 
=(\mu_1,\mu_2,\ldots ,\mu_N)$, where $\mu_i$ is the number of $i$'s 
occurring in the word $w$. We assume that the weight $\mu$ of a word $w$ 
is dominant, in other words, $\mu_1\ge\mu_2\cdots\ge\mu_N$.

i)~ First, we assume that $w$ is a standard word (i.e. its weight is 
$\mu =(1^N)$). We attach an index to each element of $w$ as follows: 
the index of 1 is 
equal to 0, and if the index of $k$ is $i$ then the index of $k+1$ is 
either $i$ or $i+1$ according as it lies to the right or left of $k$. The 
charge $c(w)$ of $w$ is defined to be the sum of the indices.

ii)~ Now assume that $w$ is a word of weight $\mu$ and $\mu$ is a 
partition. We extract a standard subword from $w$ in the following way. 
Reading $w$ from the left to the right, we choose the first occurrence of 
1, then the first occurrence of 2 to the right of the 1 chosen and so on. 
If at some step there is no $s+1$ to the right of the $s$ chosen before, 
we come back to the beginning of the word. This operation extracts from 
$w$ a standard subword $w_1$. Then we erase the word $w_1$ from $w$ and 
repeat the procedure to obtain a standard subword $w_2$, etc.

The charge of $w$ is defined as the sum of the charges of the standard 
subwords obtained in this way: $c(w)=\sum c(w_i)$. We note that the 
charge of a word $w$ is zero if and only if the word $w$ is a lattice 
word.
\begin{exm}\label{e6.1} Consider a word
\begin{equation}
w=64\underline{3}2\underline{1}1115\underline{4}33\underline{2}26
\underline{5}54\underline{6}. \label{6.6}
\end{equation}
The subword $w_1$ is 314256, consisting of the underlined symbols in $w$ 
in (\ref{6.6}). When $w_1$ is erased, we are left with
\begin{equation}
\widetilde w = 642\underline{1}11\underline{5}\underline{3}3\underline{2}
\underline{6}5\underline{4}. \label{6.7}
\end{equation}
The subword $w_2$ is 153264, consisting of the underlined 
symbols in $\widetilde 
w$ in (\ref{6.7}). When $w_2$ is erased, we are left with
\begin{equation}
6421135, \label{6.8}
\end{equation}
so that $w_3=642135$ and $w_4=1$.

Now let us compute the charges of the standard subwords $w_1$, $w_2$, 
$w_3$, $w_4$. For the word $w_1$ the indices (attached as subscripts) are 
$3_11_04_12_05_16_1$, so that $c(w_1)=4$; $1_05_23_12_06_24_1$ for $w_2$, 
so that $c(w_2)=6$; $6_34_22_11_03_15_2$ for $w_3$, so that $c(w_3)=9$; 
$1_0$ for $w_4$, so that $c(w_4)=0$: hence $c(w_4)=4+6+9+0=19$.
\end{exm}

\subsection{Spectral decomposition of Kostka-Foulkes polynomials.}

Let $\lambda$, $\mu$ be partitions of lengths $\le m$, and $d$ be a 
composition, $l(d)\le m-1$. We define the charge of a tableau 
$T\in{\mathrm{Tab}}({\mathrm{Sh}}_d(\mu ),\lambda )$ as follows
\begin{equation}
c(T)=\sum_{i=1}^{m-1}(m-i)d_i:=c(d). \label{6.9}
\end{equation}
In general, the map $\theta_d$ does not preserve the charges. However, it 
is well-known (see, e.g. \cite{LLT,K}, that the map 
$\theta_d$ does preserve the charges, if $\mu =(l^m)$ is a rectangular 
partition.
\begin{thm}\label{t6.2} Let $\mu =(l^m)$ be a rectangular partition, then
\begin{equation}
K_{\lambda ,\mu}(q)=\sum_{d=(d_1,\ldots ,d_{m-1})\in{\bf Z}_{\ge 0}^{m-1}}
q^{c(d)}|{\mathrm{LR}}_0({\mathrm{Sh}}_d(\mu ),\lambda )|. \label{6.10}
\end{equation}
\end{thm}
\begin{exm}\label{e6.3} Let us consider $\lambda =(4321)$ and $\mu 
=(2^5)$. It is easy to see that $K_{\lambda ,\mu}=24$. The Bethe ansatz 
method (see \cite{K}) gives the following expression for the 
Kostka-Foulkes polynomial
\begin{eqnarray*}
K_{\lambda ,\mu}(q)&=&q^3\left[\begin{array}{c}2\\ 1\end{array}\right]
\left[\begin{array}{c} 4\\ 1\end{array}\right]
+q^5\left[\begin{array}{c} 2\\ 1\end{array}\right]
\left[\begin{array}{c} 2\\ 1\end{array}\right]
\left[\begin{array}{c} 4\\ 1\end{array}\right]\\ \\
&=&q^3(1+2q+3q^2+5q^3+5q^4+4q^5+3q^6+q^7).
\end{eqnarray*}
Now, let us apply the spectral decomposition method,  
Theorem~\ref{t6.2}. In our example there exist 17 terms in the RHS of 
(\ref{6.10}). We give a complete list of the tableaux 
(with their charges) from the 
disjoint union $\bigsqcup_{d\in{\bf Z}_{\ge 0}^{m-1}}
{\mathrm{LR}}_0({\mathrm{Sh}}_d(\mu ),\lambda )$.
\vskip 0.2cm
$$
\begin{array}{cccc}
&\!\!\!\!\!\!\!\!&\!\!\!\!\!\!\!\!1~1\\ 
&\!\!\!\!\!\!\!\!1~2\\ 
1~2\\ 2~3\\ 3~4\\ \\7
\end{array}
\begin{array}{cccc}
&\!\!\!\!\!\!\!\!&\!\!\!\!\!\!\!\!1~1\\ 
&\!\!\!\!\!\!\!\!&\!\!\!\!\!\!\!\!2~2\\ 
&\!\!\!\!\!\!\!\!&\!\!\!\!\!\!\!\!3~3\\ 
&\!\!\!\!\!\!\!\!1~4\\ 
1~2\\ \\3
\end{array}
\hskip 2cm\begin{array}{ccccc}
&\!\!\!\!\!&\!\!\!\!\!1~1\\ 
&\!\!\!\!\!\!\!\!1~1\\ 
&\!\!\!\!\!\!\!\!2~2\\ 
&\!\!\!\!\!\!\!\!3~3\\ 
2~4\\ \\9
\end{array}
\begin{array}{ccccc}
&\!\!\!\!\!\!\!\!&\!\!\!\!\!\!\!\!1~1\\ 
&\!\!\!\!\!1~2\\ 
&\!\!\!\!\!2~3\\ 
&\!\!\!\!\!3~4\\ 
1~2\\ \\6
\end{array}
\hskip 2cm\begin{array}{cccc}
&\!\!\!\!\!&\!\!\!\!\!1~1\\ 
&\!\!\!\!\!&\!\!\!\!\!2~2\\ 
1~1\\ 2~3\\ 3~4\\ \\6
\end{array}
\begin{array}{cccc}
&\!\!\!\!\!&\!\!\!\!\!1~1\\ 
&\!\!\!\!\!&\!\!\!\!\!2~2\\ 
&\!\!\!\!\!&\!\!\!\!\!3~3\\ 
1~1\\ 2~4\\ \\4
\end{array}
$$
\vskip 0.5cm
$$
\begin{array}{cccc}
&\!\!\!\!\!\!\!\!&\!\!\!\!\!\!\!\!&\!\!\!\!\!\!\!\!1~1\\ 
&\!\!\!\!\!\!\!\!&\!\!\!\!\!\!\!\!1~2\\ 
&\!\!\!\!\!\!\!\!1~3\\ 
2~2\\ 3~4\\ \\ 9
\end{array}
\begin{array}{cccc}
&\!\!\!\!\!\!\!\!&\!\!\!\!\!\!\!\!&\!\!\!\!\!\!\!\!1~1\\ 
&\!\!\!\!\!\!\!\!&\!\!\!\!\!\!\!\!1~2\\ &\!\!\!\!\!\!\!\!2~2\\ 
1~3\\ 3~4\\ \\ 9
\end{array}
\begin{array}{cccc}
&\!\!\!\!\!\!\!\!&\!\!\!\!\!\!\!\!&\!\!\!\!\!\!\!\!1~1\\
&\!\!\!\!\!\!\!\!&\!\!\!\!\!\!\!\!&\!\!\!\!\!\!\!\!2~2\\ 
&\!\!\!\!\!\!\!\!&\!\!\!\!\!\!\!\!1~3\\ 
&\!\!\!\!\!\!\!\!1~4\\ 
2~3\\ \\ 6
\end{array}
\begin{array}{cccc}
&\!\!\!\!\!\!\!\!&\!\!\!\!\!\!\!\!&\!\!\!\!\!\!\!\!1~1\\ 
&\!\!\!\!\!\!\!\!&\!\!\!\!\!\!\!\!&\!\!\!\!\!\!\!\!2~2\\ 
&\!\!\!\!\!\!\!\!&\!\!\!\!\!\!\!\!1~3\\ 
&\!\!\!\!\!\!\!\!2~3\\ 
1~4\\ \\ 6
\end{array}
\hskip 2cm\begin{array}{cccc}
&\!\!\!\!\!\!\!\!&\!\!\!\!\!\!\!\!1~1\\ 
&\!\!\!\!\!\!\!\!&\!\!\!\!\!\!\!\!2~2\\ 
&\!\!\!\!\!\!\!\!1~3\\ 
1~2\\ 3~4\\ \\ 5
\end{array}
\begin{array}{cccc}
&\!\!\!\!\!\!\!\!&\!\!\!\!\!\!\!\!1~1\\ 
&\!\!\!\!\!\!\!\!&\!\!\!\!\!\!\!\!2~2\\ 
&\!\!\!\!\!\!\!\!1~3\\ 
1~3\\ 2~4\\ \\ 5
\end{array}
$$
\vskip 0.5cm
$$
\begin{array}{ccccc}
&\!\!\!\!\!\!\!\!&\!\!\!\!\!\!\!\!&\!\!\!\!\!\!\!\!1~1\\
&\!\!\!\!\!\!\!\!&\!\!\!\!\!\!\!\!1~2\\ 
&\!\!\!\!\!\!\!\!1~2\\ 
&\!\!\!\!\!\!\!\!3~3\\ 
2~4\\ \\ 8
\end{array}
\begin{array}{ccccc}
&\!\!\!\!\!\!\!\!&\!\!\!\!\!\!\!\!&\!\!\!\!\!\!\!\!1~1\\ 
&\!\!\!\!\!\!\!\!&\!\!\!\!\!\!\!\!1~2\\ 
&\!\!\!\!\!\!\!\!1~3\\ 
&\!\!\!\!\!\!\!\!2~4\\ 
2~3\\ \\ 8
\end{array}
\begin{array}{ccccc}
&\!\!\!\!\!\!\!\!&\!\!\!\!\!\!\!\!&\!\!\!\!\!\!\!\!1~1\\ 
&\!\!\!\!\!\!\!\!&\!\!\!\!\!\!\!\!1~2\\ 
&\!\!\!\!\!\!\!\!2~2\\ 
&\!\!\!\!\!\!\!\!3~3\\ 
1~4\\ \\ 8
\end{array}
\begin{array}{ccccc}
&\!\!\!\!\!\!\!\!&\!\!\!\!\!\!\!\!&\!\!\!\!\!\!\!\!1~1\\ 
&\!\!\!\!\!\!\!\!&\!\!\!\!\!\!\!\!1~2\\ 
&\!\!\!\!\!\!\!\!&\!\!\!\!\!\!\!\!2~3\\ 
&\!\!\!\!\!\!\!\!1~4\\ 
2~3\\ \\ 7
\end{array}
\begin{array}{ccccc}
&\!\!\!\!\!\!\!\!&\!\!\!\!\!\!\!\!&\!\!\!\!\!\!\!\!1~1\\ 
&\!\!\!\!\!\!\!\!&\!\!\!\!\!\!\!\!1~2\\ 
&\!\!\!\!\!\!\!\!&\!\!\!\!\!\!\!\!2~3\\ 
&\!\!\!\!\!\!\!\!2~3\\ 
1~4\\ \\ 7
\end{array}
\begin{array}{ccccc}
&\!\!\!\!\!\!\!\!&\!\!\!\!\!\!\!\!&\!\!\!\!\!\!\!\!1~1\\ 
&\!\!\!\!\!\!\!\!&\!\!\!\!\!\!\!\!1~2\\ 
&\!\!\!\!\!\!\!\!&\!\!\!\!\!\!\!\!2~3\\ 
&\!\!\!\!\!\!\!\!2~4\\ 
1~3\\ \\ 7
\end{array}
$$
\vskip 0.5cm
$$
\begin{array}{ccccc}
&\!\!\!\!\!&\!\!\!\!\!&\!\!\!\!\!1~1\\ 
&\!\!\!\!\!\!\!\!1~1\\ 
&\!\!\!\!\!\!\!\!2~2\\ 
2~3\\ 3~4\\ \\ 10
\end{array}
\begin{array}{ccccc}
&\!\!\!\!\!\!\!\!&\!\!\!\!\!\!\!\!&\!\!\!\!\!\!\!\!1~1\\ 
&\!\!\!\!\!\!\!\!&\!\!\!\!\!\!\!\!&\!\!\!\!\!\!\!\!2~2\\ 
&\!\!\!\!\!&\!\!\!\!\!1~3\\ 
&\!\!\!\!\!&\!\!\!\!\!3~4\\ 
1~2\\ \\ 5
\end{array}
\hskip 2cm\begin{array}{cccc}
&\!\!\!\!\!\!\!\!&\!\!\!\!\!\!\!\!&\!\!\!\!\!\!\!\!1~1\\ 
&\!\!\!\!\!&\!\!\!\!\!1~2\\ 
&\!\!\!\!\!&\!\!\!\!\!2~3\\ 
1~2\\ 3~4\\ \\ 8
\end{array}
\begin{array}{cccc}
&\!\!\!\!\!&\!\!\!\!\!&\!\!\!\!\!1~1\\ 
&\!\!\!\!\!&\!\!\!\!\!&\!\!\!\!\!2~2\\ 
&\!\!\!\!\!\!\!1~1\\ 
&\!\!\!\!\!\!\!\!3~3\\ 
2~4\\ \\ 7
\end{array}
\hskip 2cm\begin{array}{cccc}
&\!\!\!\!\!\!\!\!&\!\!\!\!\!\!\!\!1~1\\ 
&\!\!\!\!\!\!\!\!1~2\\ 
&\!\!\!\!\!\!\!\!2~3\\ 
1~3\\ 2~4\\ \\ 6
\end{array}
\begin{array}{cccc}
&\!\!\!\!\!\!\!\!&\!\!\!\!\!\!\!\!1~1\\ 
&\!\!\!\!\!\!\!\!&\!\!\!\!\!\!\!\!2~2\\ 
&\!\!\!\!\!\!\!\!1~3\\ 
&\!\!\!\!\!\!\!\!2~4\\ 
1~3\\ \\ 4 \\
\end{array}
$$
Hence, the spectral decomposition of the Kostka-Foulkes polynomial 
$K_{(4321),(2^5)}(q)$ has the following form
$$q^3+(1+1)q^4+(2+1)q^5+(1+2+1+1)q^6+(1+3+1)q^7+(1+3)q^8+(1+2)q^9+q^{10}.
$$
One can recognize a certain symmetry among the tableaux from the above
list. This symmetry is a consequence of the Sch\"utzenberger involution 
action on the set ${\rm SST}(\lambda ,\mu )$. More precisely, it follows 
from \cite{K}, Corollary~4.2, that if we denote by $S$ the 
Sch\"utzenberger involution
$$S:~{\rm SST}(\lambda ,(l^m))\to{\rm SST}(\lambda ,(l^m)),
$$
then we have the following equality for the exponents:
$$d(S(T))={\overleftarrow d}(T),
$$
where for any sequence $\alpha =(\alpha_1,\ldots ,\alpha_m)$ we set 
${\overleftarrow\alpha}:=(\alpha_m,\alpha_{m-1},\ldots ,\alpha_2,\alpha_1)$.
On the other hand, it follows from Corollary~\ref{c5.3} that
$$|{\mathrm{LR}}_0({\mathrm{Sh}}_d(\mu ),\lambda )|=
|{\mathrm{LR}}_0({\mathrm{Sh}}_{\overleftarrow d}(\mu ),\lambda )|.
$$
If we define ${\mathrm{ind}}(d)=\sum_{j=1}^{m-1}jd_j$
(cf.\ \cite{K,KR}), we obtain the 
following
\begin{cor}\label{c6.3} Let $\mu =(l^m)$ be a rectangular partition, then
$$K_{\lambda ,\mu}(q)=\sum_{d=(d_1,\ldots ,d_{m-1})\in{\bf Z}_{\ge 
0}^{m-1}} q^{{\mathrm{ind}}(d)}
|{\mathrm{LR}}_0({\mathrm{Sh}}_d(\mu ),\lambda )|.
$$
\end{cor}
\end{exm}

\section{Truncated characters and
branching functions}\label{sec:truncated}

\subsection{Kostka--Foulkes polynomials and the
truncated characters}\label{final}

We shall explain how our main theorems,
Theorems \ref{thm:maintheorem} and \ref{t6.2}, 
are related to each other.
Let us consider the special situation
of  (\ref{eq:djkmo}) where
$\Lambda(K)=l\Lambda_k$, $(k=0,\dots,n-1)$.
The configuration space ${\cal S}_k:=
{\cal S}_K$ admits
a filtration,
${\cal S}_{k,k}\subset
{\cal S}_{k,n+k}\subset
{\cal S}_{k,2n+k}\subset
\cdots
\subset
{\cal S}_k$,
where
$$
{\cal S}_{k,m}=\{
(s_1,\dots,s_{m},(v_{n\ldots n},
v_{n-1\ldots n-1},\dots,
v_{1\ldots 1})^\infty)\},
\qquad
m \equiv k\ {\mathrm{mod}}\ n.
$$
Each element 
 of ${\cal S}_{k,m}$ is
naturally identified
with a spin configuration 
on a finite lattice of size $m$ by the truncation
$(s_i)_{i=1}^\infty\mapsto (s_i)_{i=1}^m $.
For a tableau $T$ we define the $\hbox{\germ gl}_n$-weight 
of $T$ as 
$\overline{\mathrm{wt}}(T):= \sum_{a=1}^n m_a \overline{\epsilon}_a$,
where $(m_1, \ldots, m_n)$ is the content of $T$ and 
$\overline{\epsilon}_1, \ldots, \overline{\epsilon}_n$ are 
linearly independent vectors.
We introduce the {\it truncated character} $F_{k,m}$ of ${\cal S}_k$,
\begin{align*}
F_{k,m}(q,x) &=q^{A_{k,m}}
\sum_{\vec{s}\in {\cal S}_{k,m}}
q^{E(\vec{s})}
e^{\sum_{i=1}^m {\overline{\mathrm{wt}}(s_i)}},\cr
A_{k,m}&={1\over2}lnN(N-1)+lNk,
\quad
N={1\over n}(m-k),
\end{align*}
where $x_1 = e^{\overline{\epsilon}_1}, \ldots, 
x_n = e^{\overline{\epsilon}_n}$ are
independent variables here.
By definition, it has the property,
\begin{equation}\label{chf}
{\mathrm{ch}}\, {\cal L}(l \Lambda_k)(q,x) 
=\lim_{{m\rightarrow \infty\atop
m\equiv k\, {\mathrm{mod}}\, n}}
q^{-A_{k,m}}
\left. F_{k,m}(q,x) \right|_{x_1 \cdots x_n = 1}.
\end{equation}
Here we used the fact that $H_l(v_{n\ldots n},s)(=l)$ is
independent of $s\in B_l$.
The following theorem states that the branching
functions of the truncated character
are the Kostka--Foulkes
polynomials.
\begin{thm}[\cite{NY}]
\label{thm:finite}
\begin{equation}
\label{eq:finitedecom}
F_{k,m}(q,x)
=
\sum_{\lambda} K_{\lambda,(l^m)}(q)
 s_\lambda(x),
\qquad
m\equiv k\ {\mathrm{mod}}\ n,
\end{equation}
where the sum extends over the partitions $\lambda$ of $lm$ 
such that $l(\lambda) \le n$.
\end{thm}

The spectral decomposition 
of ${\mathrm{ch}}\, {\cal L}(l\Lambda_k)$
(Theorem \ref{thm:maintheorem})
induces that of $F_{k,m}$,
\begin{equation}
\label{eq:intde}
F_{k,m}(q,x) 
=
\sum_{d=(d_1,\dots,d_{m-1})}
q^{\sum_{i=1}^{m-1}
(m-i) d_i}
\overline{t}_{{\mathrm{Sh}}_d((l^m))}(x),\quad
\overline{t}_{\lambda/\mu}(x):= \sum_{T \in {\mathrm{NMT}}(\lambda/\mu)}
e^{\overline{\mathrm{wt}}(T)}.
\end{equation}
Furthermore,
the subcharacter $\overline{t}_{{\mathrm{Sh}}_d((l^m))}$ is expanded 
by the
Schur functions as 
\begin{equation}\label{eq:subd}
\overline{t}_{{\mathrm{Sh}}_d((l^m))}
=
\sum_{\lambda}
| {\mathrm{LR}}_0({\mathrm{Sh}}_d((l^m)),\lambda)|
s_\lambda
\end{equation}
thanks to  Propositions \ref{prop:characterformula} and 
 \ref{p5.2}.
By comparing 
the formulae (\ref{eq:finitedecom})--(\ref{eq:subd}),
we see that
 Theorem \ref{t6.2}
is equivalent to
 Theorem \ref{thm:finite}.\footnote{Thus we
 have given an alternative proof 
of Theorem \ref{thm:finite}.
Conversely, we
have our second proof of Theorem \ref{t6.2} using
Theorem \ref{thm:finite}.}
In other words, we can interpret
Theorem \ref{t6.2} 
as the decomposition of
the branching functions $K_{\lambda,(l^m)}(q)$
induced from the spectral decomposition of
the truncated character $F_{k,m}$.

Below we remark the relation between our bijection,
\begin{equation}
\theta : {\mathrm{SST}}(\lambda,(l^m))
\longrightarrow \bigsqcup_d {\mathrm{LR}}_0
({\mathrm{Sh}}_d((l^m)),\lambda),\qquad
\theta|_{{\mathrm{SST}}_d(\lambda,(l^m))}=\theta_d,
\end{equation}
and the bijection by Nakayashiki and Yamada  \cite{NY},
\begin{align*}
&\pi : {\mathrm{SST}}(\lambda,(l^m))
\longrightarrow (B_l^m)_\lambda^{\mathrm{high}},\\
&(B_l^m)_\lambda^{\mathrm{high}}:=
\{
\vec{s}=(s_1,\dots,s_m)\in B_l^{m} \mid
\sum_{i=1}^m {\mathrm{wt}}(s_i)=\lambda,\quad
\tilde{e}_i(s_1\otimes \dots \otimes s_m)=0
\}.
\end{align*}
Here $\pi(T)=(s_1,\dots,s_m)$ is defined so that
the contents of $s_i$ are the indices
of the rows in $T$ which contain
$i$.
For example,
\begin{equation*}
\begin{picture}(300,30)(0,10)
\put(30,10){\line(1,0){20}}
\put(30,20){\line(1,0){30}}
\put(30,30){\line(1,0){40}}
\put(30,40){\line(1,0){40}}
\put(125,20){\line(1,0){30}}
\put(125,30){\line(1,0){30}}
\put(193,20){\line(1,0){30}}
\put(193,30){\line(1,0){30}}
\put(261,20){\line(1,0){30}}
\put(261,30){\line(1,0){30}}
\put(30,10){\line(0,1){30}}
\put(40,10){\line(0,1){30}}
\put(50,10){\line(0,1){30}}
\put(60,20){\line(0,1){20}}
\put(70,30){\line(0,1){10}}
\put(125,20){\line(0,1){10}}
\put(135,20){\line(0,1){10}}
\put(145,20){\line(0,1){10}}
\put(155,20){\line(0,1){10}}
\put(193,20){\line(0,1){10}}
\put(203,20){\line(0,1){10}}
\put(213,20){\line(0,1){10}}
\put(223,20){\line(0,1){10}}
\put(261,20){\line(0,1){10}}
\put(271,20){\line(0,1){10}}
\put(281,20){\line(0,1){10}}
\put(291,20){\line(0,1){10}}
\put(3,23){$T = $}
\put(75,23){$\mapsto \bigl( s_1 = $}
\put(157,18){$,$}
\put(163,23){$s_2 = $}
\put(225,18){$,$}
\put(231,23){$s_3 = $}
\put(297,23){$\bigr).$}
\put(32,31){\small 1}
\put(32,21){\small 2}
\put(32,11){\small 3}
\put(42,31){\small 1}
\put(42,21){\small 2}
\put(42,11){\small 3}
\put(52,31){\small 1}
\put(52,21){\small 2}
\put(62,31){\small 3}
\put(127,21){\small 1}
\put(137,21){\small 1}
\put(147,21){\small 1}
\put(195,21){\small 2}
\put(205,21){\small 2}
\put(215,21){\small 2}
\put(263,21){\small 1}
\put(273,21){\small 3}
\put(283,21){\small 3}
\end{picture}
\end{equation*}
Let
\begin{align*}
(B_l^m)_\lambda &:=
\{\vec{s}=(s_1,\dots,s_m)\in B_l^{m} \mid
\sum_{i=1}^m {\mathrm{wt}}(s_i)=\lambda \},\\
{\mathrm{NMT}}
({\mathrm{Sh}}_d((l^m)),\lambda) &:= 
\{ T \in {\mathrm{NMT}}(
{\mathrm{Sh}}_d((l^m))) \mid 
\mbox{the content of $T$ is $\lambda$} \}.
\end{align*}
We have a natural bijection
\begin{equation}
\Phi : (B_l^m)_\lambda
\longrightarrow
\bigsqcup_d {\mathrm{NMT}}
({\mathrm{Sh}}_d((l^m)),\lambda)
\end{equation}
such that the content of the $i$th row
of $\Phi(s_1,\dots,s_m)$ equals to that of
$s_i$;
for example,
\begin{equation*}
\begin{picture}(235,30)(5,10)
\put(40,20){\line(1,0){30}}
\put(40,30){\line(1,0){30}}
\put(80,20){\line(1,0){30}}
\put(80,30){\line(1,0){30}}
\put(120,20){\line(1,0){30}}
\put(120,30){\line(1,0){30}}
\put(195,10){\line(1,0){30}}
\put(195,20){\line(1,0){40}}
\put(205,30){\line(1,0){30}}
\put(205,40){\line(1,0){30}}
\put(40,20){\line(0,1){10}}
\put(50,20){\line(0,1){10}}
\put(60,20){\line(0,1){10}}
\put(70,20){\line(0,1){10}}
\put(80,20){\line(0,1){10}}
\put(90,20){\line(0,1){10}}
\put(100,20){\line(0,1){10}}
\put(110,20){\line(0,1){10}}
\put(120,20){\line(0,1){10}}
\put(130,20){\line(0,1){10}}
\put(140,20){\line(0,1){10}}
\put(150,20){\line(0,1){10}}
\put(195,10){\line(0,1){10}}
\put(205,10){\line(0,1){30}}
\put(215,10){\line(0,1){30}}
\put(225,10){\line(0,1){30}}
\put(235,20){\line(0,1){20}}
\put(5,23){$\Phi : \bigl( $}
\put(72,18){$,$} 
\put(112,18){$,$} 
\put(160,23){$\bigr)\, \mapsto $}
\put(42,21){\small 1}
\put(52,21){\small 1}
\put(62,21){\small 1}
\put(82,21){\small 2}
\put(92,21){\small 2}
\put(102,21){\small 2}
\put(122,21){\small 1}
\put(132,21){\small 3}
\put(142,21){\small 3}
\put(197,11){\small 1}
\put(207,11){\small 3}
\put(207,21){\small 2}
\put(207,31){\small 1}
\put(217,11){\small 3}
\put(217,21){\small 2}
\put(217,31){\small 1}
\put(227,31){\small 1}
\put(227,21){\small 2}
\put(239,19){$.$}
\end{picture}
\end{equation*}
Then

\begin{prop}
The diagram
\begin{equation*}
\begin{matrix}
 (B^m_l)^{\rm high}_\lambda & \subset & (B^m_l)_\lambda \\
\qquad\overset{\pi\  \simeq}{} {\nearrow}\hfill & & \\
{\mathrm{SST}}(\lambda, (l^m))\hfill & &  
\downarrow \overset{\simeq \ \Phi}{} \\
\qquad\overset{\theta \ \simeq}{} {\searrow} \hfill & & \\
\qquad\bigsqcup_d {\mathrm{LR}}_0
({\mathrm{Sh}}_d((l^m)),\lambda) & \subset &
\bigsqcup_d {\mathrm{NMT}}
({\mathrm{Sh}}_d((l^m)),\lambda) \\
\end{matrix}
\end{equation*}
is commutative.
\end{prop}

\begin{cor}
$\Phi$ gives a bijection between
$(B_l^m)_\lambda^{\mathrm{high}}$
and\/ $\bigsqcup_d {\mathrm{LR}}_0
({\mathrm{Sh}}_d((l^m)),\lambda)$.
\end{cor}

\subsection{Branching functions as limits of the Kostka--Foulkes
polynomials}

As an application of Theorem \ref{thm:finite} we can 
identify a limit of the Kostka--Foulkes polynomial 
$K_{\lambda, (l^m)}(q)$ with a branching function
$b^\Lambda_\lambda(q)$ of an $\widehat{sl}_n$-module:
\begin{equation}
{\rm ch}\, {\cal L}(\Lambda )(q,x)=
\sum_{\lambda}b_{\lambda }^{\Lambda}(q)
\left. s_{\lambda} (x) \right|_{x_1 \cdots x_n = 1}, 
\label{6.11}
\end{equation}
where the sum extends over all partitions 
$\lambda$ such that $l(\lambda )\le n-1$, and 
$|\lambda |\equiv |\Lambda | \mod n$.
By comparing (\ref{chf}), (\ref{eq:finitedecom}) and 
(\ref{6.11}), we have 
\begin{cor}[\cite{NY},\cite{K2}]\label{c6.2}
Let $0\le k\le n-1$. Then
\begin{equation}
b_{\lambda}^{l\Lambda_k}(q)=
\lim_{N\to\infty}q^{-A_{k,k+Nn}}K_{\lambda_N,\mu_N}(q), 
\label{6.12}
\end{equation}
where $\lambda_N=\lambda +
\left(\left( \frac{
l(k+Nn)-|\lambda 
|}{n}\right)^n\right)$
and $\mu_N=(l^{k+Nn})$.
\end{cor}
The formula (\ref{6.12}) was originally conjectured in \cite{K2}
for $k=0$ and proved in \cite{NY}.
A similar formula to Corollary \ref{c6.2}
for the branching functions $b^\Lambda_\lambda$ 
for a general
level $l$ dominant integral weight
$\Lambda\neq l\Lambda_k$
can be obtained, as we describe below, by
considering the truncated character of the
space ${\cal S}_K$, $\Lambda=\Lambda(K)$.
In the general case, 
however, the polynomials appearing
in the right hand side of (\ref{6.12})
are no longer identified with the Kostka--Foulkes
polynomials.

For a given semistandard tableau of
skew shape $T\in {\mathrm{SST}}(\lambda/\nu,\mu)$,
$l(\mu)\leq m$,
we define the (extended) exponents
$d_0(T), d_1(T),\dots,d_{m-1}(T)$,
where the definition of $d_i(T)$, $i>0$ remains
the same as in Definition \ref{d5.1},
and 
$d_0 (T):= \mu_1 - \zeta_0(T)$ 
with $\zeta_0(T)$ the number of $1$'s
in $T$ which does {\it not\/} belong to the
first row.
Notice that, when $\lambda/\nu$ is a Young
diagram, i.e., $\nu=\emptyset$,
$d_0(T)=\mu_1$ is independent of $T$ for
given $\lambda$ and $\mu$.
Let ${\mathrm{SST}}_d(\lambda/\nu 
,\mu )$ be the
set of all semistandard tableaux $T\in
{\mathrm{SST}}(\lambda/\nu ,\mu )$
such that $d_i(T)=d_i$, $0\le i\le m-1$. 
Theorem \ref{t5.1} is
extended to skew shape tableaux as
\begin{thm}\label{thm:extbij}
There is a natural bijection
\begin{equation}
\theta_d: {\mathrm{SST}}_d(\lambda/\nu ,\mu )
\leftrightarrow
{\mathrm{LR}}_0({\mathrm{Sh}}_d(\mu,\nu ),\lambda ),
\end{equation}
where $d=(d_0,d_1,\dots,d_{m-1})$,
\begin{equation*}
\begin{picture}(224,105)(-62,-15)
\put(0,0){\line(1,0){60}}
\put(0,10){\line(1,0){80}}
\put(20,20){\line(1,0){60}}
\put(60,40){\line(1,0){60}}
\put(60,50){\line(1,0){80}}
\put(80,60){\line(1,0){80}}
\put(100,70){\line(1,0){20}}
\put(120,80){\line(1,0){20}}
\put(140,90){\line(1,0){20}}
\put(0,10){\line(0,-1){10}}
\put(60,10){\line(0,-1){10}}
\put(20,20){\line(0,-1){10}}
\put(80,20){\line(0,-1){10}}
\put(60,50){\line(0,-1){10}}
\put(120,50){\line(0,-1){10}}
\put(80,60){\line(0,-1){10}}
\put(140,60){\line(0,-1){10}}
\put(100,70){\line(0,-1){10}}
\put(120,80){\line(0,-1){10}}
\put(140,90){\line(0,-1){10}}
\put(160,90){\line(0,-1){30}}
\put(0,10){\small$\overbrace{\hskip0pt}$}
\put(-5,20){\small$ d_{m-1}$}
\put(60,50){\small$\overbrace{\hskip0pt}$}
\put(65,60){\small$ d_{1}$}
\put(80,60){\small$\overbrace{\hskip0pt}$}
\put(85,70){\small$ d_{0}$}
\put(0,0){\small$\underbrace{\hskip60pt}$}
\put(27,-15){\small$ \mu_{m}$}
\put(20,20){\small$\overbrace{\hskip60pt}$}
\put(42,30){\small$ \mu_{m-1}$}
\put(80,35){\circle*{2}}
\put(75,30){\circle*{2}}
\put(140,70){$\widehat \nu$}
\put(-62,50){${\mathrm{Sh}}_d(\mu,\nu )=$}
\put(162,50){,}
\end{picture}
\end{equation*}
and $\widehat \nu$ is the skew diagram
obtained from $\nu$ by the 180 degree rotation.
\end{thm}
The bijection $\theta_d$ is defined as follows:
For a given $T\in {\mathrm{SST}}_d(\lambda/\nu ,\mu )$,
the contents of the $l(\nu)+k$th row 
in $\theta_d(T)$ are given by the indices of the
rows in $T$ which contain $k$
(the contents of the first $l(\nu)$ rows in $\theta_d(T)$
are uniquely determined from the semistandardness and the
LR property of $\theta_d(T)$). 
The theorem is easily derived from
 Theorem \ref{t5.1}, and it includes
Theorem \ref{t5.1} as the special case
when $\nu=\emptyset$.

Keeping Theorem \ref{thm:extbij} in mind,
we introduce two functions.
Firstly, for a skew diagram $\lambda/\nu$
and a composition $\mu$ with
$|\lambda/\nu|=|\mu|$, $m=l(\mu)$,
we define a polynomial,
\begin{equation*}
G_{\lambda/\nu,\, \mu}(q):=
\sum_{T\in {\mathrm{SST}}(\lambda/\nu,\, \mu)}
q^{\sum_{i=0}^{m-1}(m-i)d_i(T)}.
\end{equation*}
Secondly, for 
$K=(k_i)\in {\cal K}_l$ and a positive
integer $m$, we define the truncated
character $F_{K,m}$ of ${\cal S}_K$
as follows: Consider a filtration
${\cal S}_{K,1}\subset
{\cal S}_{K,2}\subset
\cdots
\subset
{\cal S}_{K}$, where
\begin{equation*}
{\cal S}_{K,m}:=
\{ (s_i) \mid
H(s_{i+1},s_i)=k_i\ \mbox{for any
$i\geq m+1$}\}.
\end{equation*}
For $\vec{s}\in {\cal S}_{K,m}$,
$s_i=s_i^{(K)}$ holds for any $i\geq m+n$,
but not necessarily so for $i \leq m+n-1$
(cf.\ Fig.\ \ref{fig:finite}).
Now we define
\begin{align*}
F_{K,m}(q,x):=&
q^{B_{K,m}}\left(
\prod_{i=2}^n x_i^{-w_i}
\right)
\sum_{\vec{s}\in {\cal S}_{K,m}}
q^{E(\vec{s})}
e^{\sum_{i=1}^{m+n-1} \overline{{\mathrm{wt}}}(s_i)},\cr
B_{K,m}=&\sum_{i=1}^m i k_i,
\qquad
w_i=\sum_{j=1}^i k_{m+j-1}.
\end{align*}
Then, due to the spectral decomposition of 
${\cal S}_K$ and Theorem \ref{t5.1},
we have a generalization of
Theorem \ref{thm:finite}.
\begin{thm}
\begin{equation*}
F_{K,m}(q,x)
=
\sum_{\lambda} G_{\lambda/\nu_{K,m},(l^m)}(q) s_\lambda(x),
\end{equation*}
where
\begin{equation*}
\begin{picture}(127,55)(-45,-15)
\put(0,0){\line(1,0){20}}
\put(20,10){\line(1,0){20}}
\put(40,20){\line(1,0){10}}
\put(50,30){\line(1,0){20}}
\put(0,40){\line(1,0){70}}
\put(0,40){\line(0,-1){40}}
\put(20,10){\line(0,-1){10}}
\put(40,20){\line(0,-1){10}}
\put(50,30){\line(0,-1){10}}
\put(70,40){\line(0,-1){10}}
\put(1,0){\small$\underbrace{\hskip0pt}$}
\put(3,-15){\small$ k_{m+n-1}$}
\put(21,10){\small$\underbrace{\hskip0pt}$}
\put(23,-5){\small$ k_{m+n-2}$}
\put(51,30){\small$\underbrace{\hskip0pt}$}
\put(53,15){\small$  k_{m+1}$}
\put(-45,15){$\nu_{K,m}=$}
\put(82,15){,}
\end{picture}
\end{equation*}
and the sum extends over the partitions $\lambda$ of
 $lm+|\nu_{K,m}|$ 
such that $l(\lambda) \le n$.
\end{thm}
Notice that $\nu_{K,m}=\emptyset$
if and only if $\Lambda(K)=l\Lambda_k$,
and $m \equiv k \mod n$ for some $k$,
$0\leq k \leq n-1$, where
the theorem reduces to Theorem \ref{thm:finite}.

\begin{cor}\label{cor:branchtwo}
Let $b^{\Lambda}_\lambda(q)$ be the branching
function of the level $l$ integrable module
${\cal L}(\Lambda)$, $\Lambda=\Lambda(K)$.
Then
\begin{equation*}
b^{\Lambda}_\lambda(q)=
\lim_{m \to \infty} q^{-B_{K,m}}
G_{\lambda_m/\nu_{K,m},(l^m)}(q),
\end{equation*}
where $B_{K,m}=\sum_{i=1}^m i k_i$ and
$\lambda_m=\lambda+\left(
\left({lm+|\nu_{K,m}|-|\lambda|\over n}\right)^n\right)$.

\end{cor}

\end{document}